\documentclass[aps,pra,letterpaper,notitlepage,longbibliography,reprint,onecolumn,superscriptaddress]{revtex4-1}
\usepackage[table]{xcolor}
\usepackage{graphicx,amsmath,amsfonts,amssymb,amsthm,bm,bbm}
\usepackage{tikz,pgfplots}
\pgfplotsset{compat=newest}
\pgfplotsset{plot coordinates/math parser=false}
\usetikzlibrary{arrows,decorations.pathmorphing,backgrounds,positioning,fit,calc,patterns,shapes,decorations.markings,shadings,snakes}
\usetikzlibrary{external}
\usepackage{ifthen}
\usepackage{multirow}
\usepackage{diagbox}
\usepackage{stmaryrd}
\usepackage{tabularx}
\usepackage{booktabs}

\newcommand{\greycline}[1]{\arrayrulecolor[gray]{.8}\cline{#1}\arrayrulecolor{black}}

\newtheoremstyle{exampleA}
{\topsep}
{10pt}
{}
{}
{\bfseries}
{\\\\}
{ 10mm }
{\thmname{#1}\thmnumber{ #2}\thmnote{ (#3)}}

\theoremstyle{exampleA}
\newtheorem{exmp}{Example}[section]

\newtheoremstyle{questionA}
{1pt}
{5pt}
{}
{}
{\bfseries}
{: }
{ 5mm }
{\thmname{#1} \thmnumber{#2}}
\theoremstyle{questionA}
\newtheorem{question}{Question}

\usepackage[pdfpagelabels,pdftex,bookmarks,breaklinks]{hyperref}
\definecolor{darkblue}{RGB}{0,0,127} 
\definecolor{darkgreen}{RGB}{0,180,0}
\hypersetup{
	colorlinks, 
	linkcolor=darkblue, 
	citecolor=darkgreen, 
	filecolor=red, 
	urlcolor=blue,
	pdftitle={Domain walls in topological phases and the Brauer-Picard ring for Vec(Zp)}, 
	pdfauthor={Daniel Barter, Jacob C Bridgeman, and Corey Jones}
}

\definecolor{tctwistcolor}{RGB}{255,0,0}
\definecolor{tcmcolor}{RGB}{0,0,180}
\definecolor{tcecolor}{RGB}{0,180,0}
\definecolor{tcTTppcolor}{RGB}{184,134,11}
\definecolor{tcTTpmcolor}{RGB}{252,15,192}

\theoremstyle{plain}
\newtheorem{theorem}{Theorem}

\newtheorem{proposition}[theorem]{Proposition}

\theoremstyle{definition}
\newtheorem{definition}[theorem]{Definition}

\newcommand{\ZZ}[1]{\mathbb{Z}/#1\mathbb{Z}}

\let\oldonlinecite\onlinecite
\renewcommand{\onlinecite}[1]{Ref.~[\oldonlinecite{#1}]}

\newcommand{\D}[1]{\mathcal{D}{} }

\newcommand{\restrict}[1]{\raise-.2ex\hbox{\ensuremath|}_{#1}}

\definecolor{tensorblue}{rgb}{0.8,0.8,1}
\definecolor{tensorred}{rgb}{1,0.5,0.5}
\definecolor{tensorgreen}{rgb}{0.6,1,0.6}
\definecolor{tensorpurp}{rgb}{1,0.5,1}

\tikzset{ten/.style={fill=tensorblue}}
\tikzset{tenred/.style={fill=tensorred}}
\tikzset{tengreen/.style={fill=tensorgreen}}
\tikzset{tenpurp/.style={fill=tensorpurp}}

\makeatletter
\newcommand{\vast}{\bBigg@{4}}
\newcommand{\Vast}{\bBigg@{9}}
\makeatother

\def\Put(#1,#2)#3{\leavevmode\makebox(0,0){\put(#1,#2){#3}}}

\makeatletter
\def\pgf@plot@curveto@handler@finish{%
  \ifpgf@plot@started%
    \pgfpathcurvebetweentimecontinue{0}{0.995}{\pgf@plot@curveto@first}{\pgf@plot@curveto@first@support}{\pgf@plot@curveto@second}{\pgf@plot@curveto@second}%
  \fi%
}
\makeatother

\newlength\figureheight 
\newlength\figurewidth 
  
\newcommand{\includeTikz}[2]
{
	\tikzifexternalizing
	{
		\includeTikzrm{#1}{#2}
	}
	{
		\IfFileExists{figures/#1.pdf}{
			\includegraphics{figures/#1}
		}
		{
			\includeTikzrm{#1}{#2}
		}
	}
}

\newcommand{\includeTikzrm}[2]{
	\tikzset{external/remake next}
	\tikzsetnextfilename{#1}
	#2
}
  
  \definecolor{nicegreena}{RGB}{1,115,16}
  \definecolor{nicegreenb}{RGB}{1,240,16}
  
  \definecolor{nicegreen}{RGB}{60,183,82}
  
    \colorlet{ccred}{red!20}
    \colorlet{ccgreen}{green!50}
    \colorlet{ccblue}{blue!20}

  \tikzset{hexr/.style= {shape=regular polygon,regular polygon sides=6,minimum size=1cm, draw,inner sep=0,anchor=center,fill=red!50}}
  \tikzset{hexg/.style= {shape=regular polygon,regular polygon sides=6,minimum size=1cm, draw,inner sep=0,anchor=center,fill=green!50}}
  \tikzset{hexb/.style= {shape=regular polygon,regular polygon sides=6,minimum size=1cm, draw,inner sep=0,anchor=center,fill=blue!50}}
  
  \definecolor{tensorblue}{rgb}{0.8,0.8,1}
  \definecolor{tensorred}{rgb}{1,0.5,0.5}
  \definecolor{tensorpurp}{rgb}{1,0.5,1}
  
  \tikzset{nonesty/.style={fill=none,draw=none}}
  \tikzset{ten/.style={fill=tensorblue}}
  \tikzset{tenred/.style={fill=tensorred}}
  \tikzset{tengreen/.style={fill=green!50!black!50}}
  \tikzset{tenpurp/.style={fill=tensorpurp}}
  \tikzset{tengrey/.style={fill=black!20}}
  \tikzset{tenorange/.style={fill=orange!30}}
  \tikzset{u/.style={fill=blue!20,draw=black}}
  \tikzset{w/.style={fill=green!50!black!50,draw=black}}

\tikzset{external/system call={pdflatex \tikzexternalcheckshellescape -halt-on-error -interaction=batchmode -jobname "\image" "\texsource"; ps2pdf13 -dEmbedAllFonts=true -r100000  "\image".pdf "\image-13".pdf && cp "\image-13".pdf "\image".pdf && rm "\image-13".pdf && rm "\image".log && rm "\image".dpth && rm "\image"Notes.bib}}

\tikzstyle{inline text}=[text height=1.1ex, text depth=0.1ex, yshift=-.1ex]

\def\ladder[#1][#2][#3][#4][#5][#6]{
	\path(#1);
	\pgfgetlastxy{\XCoord}{\YCoord};
	\begin{scope}[xshift=\XCoord,yshift=\YCoord]
		\draw (-.6,-.5)--(-.6,.5);
		\draw (.6,-.5)--(.6,.5);
		\ifthenelse{\equal{#4}{}}{}
			{
				\draw (-.6,-.25)--(.6,.25);
			};
		\node[below,inline text] at (-.6,-.5) {\footnotesize#2};
		\node[below,inline text] at (.6,-.5) {\footnotesize#3};
		\node[below,inline text] at (0,0) {\footnotesize#4};
		\node[above,inline text] at (-.6,.5) {\footnotesize#5};
		\node[above,inline text] at (.6,.5) {\footnotesize#6};
	\end{scope}
}

\def\Lladder[#1][#2][#3][#4][#5][#6][#7]{
	\path(#1);
	\pgfgetlastxy{\XCoord}{\YCoord};
	\begin{scope}[xshift=\XCoord,yshift=\YCoord]
		\draw (-1.8,-.5)--(-1.8,.5);
		\draw (-.6,-.5)--(-.6,.5);
		\draw (.6,-.5)--(.6,.5);
		\ifthenelse{\equal{#4}{}}{}
		{
			\draw (-.6,-.25)--(.6,.25);
		};
		\node[below,inline text] at (-.6,-.5) {\footnotesize#2};
		\node[below,inline text] at (.6,-.5) {\footnotesize#3};
		\node[below,inline text] at (0,0) {\footnotesize#4};
		\node[above,inline text] at (-.6,.5) {\footnotesize#5};
		\node[above,inline text] at (.6,.5) {\footnotesize#6};
		\node[below,inline text] at (-1.8,-.5) {\footnotesize#7};
	\end{scope}
}

\def\Lmorphism[#1][#2][#3][#4][#5][#6][#7]{
	\path(#1);
	\pgfgetlastxy{\XCoord}{\YCoord};
	\begin{scope}[xshift=\XCoord,yshift=\YCoord]
		\draw (-1.8,-.5)to[out=90,in=190](-.6,.25);
		\draw (-.6,-.5)--(-.6,.5);
		\draw (.6,-.5)--(.6,.5);
		\ifthenelse{\equal{#4}{}}{}
		{
			\draw (-.6,-.25)--(.6,.25);
		};
		\node[below,inline text] at (-.6,-.5) {\footnotesize#2};
		\node[below,inline text] at (.6,-.5) {\footnotesize#3};
		\node[below,inline text] at (0,0) {\footnotesize#4};
		\node[above,inline text] at (-.6,.5) {\footnotesize#5};
		\node[above,inline text] at (.6,.5) {\footnotesize#6};
		\node[below,inline text] at (-1.8,-.5) {\footnotesize#7};
	\end{scope}
}

\def\Rladder[#1][#2][#3][#4][#5][#6][#7]{
	\path(#1);
	\pgfgetlastxy{\XCoord}{\YCoord};
	\begin{scope}[xshift=\XCoord,yshift=\YCoord]
		\draw (1.8,-.5)--(1.8,.5);
		\draw (-.6,-.5)--(-.6,.5);
		\draw (.6,-.5)--(.6,.5);
		\ifthenelse{\equal{#4}{}}{}
		{
			\draw (-.6,-.25)--(.6,.25);
		};
		\node[below,inline text] at (-.6,-.5) {\footnotesize#2};
		\node[below,inline text] at (.6,-.5) {\footnotesize#3};
		\node[below,inline text] at (0,0) {\footnotesize#4};
		\node[above,inline text] at (-.6,.5) {\footnotesize#5};
		\node[above,inline text] at (.6,.5) {\footnotesize#6};
		\node[below,inline text] at (1.8,-.5) {\footnotesize#7};
	\end{scope}
}

\def\Rmorphism[#1][#2][#3][#4][#5][#6][#7]{
	\path(#1);
	\pgfgetlastxy{\XCoord}{\YCoord};
	\begin{scope}[xshift=\XCoord,yshift=\YCoord]
		\draw (1.8,-.5)to[out=90,in=-10](.6,.3);
		\draw (-.6,-.5)--(-.6,.5);
		\draw (.6,-.5)--(.6,.5);
		\ifthenelse{\equal{#4}{}}{}
		{
			\draw (-.6,-.35)--(.6,.15);
		};
		\node[below,inline text] at (-.6,-.5) {\footnotesize#2};
		\node[below,inline text] at (.6,-.5) {\footnotesize#3};
		\node[below,inline text] at (0,-.1) {\footnotesize#4};
		\node[above,inline text] at (-.6,.5) {\footnotesize#5};
		\node[above,inline text] at (.6,.5) {\footnotesize#6};
		\node[below,inline text] at (1.8,-.5) {\footnotesize#7};
	\end{scope}
}

\def\Laction[#1][#2][#3]{
		\draw (0,-.5)--(0,.5);
		\draw (-.6,-.5)to[out=90,in=210](0,0);
		\node[below,inline text] at (-.6,-.5) {\footnotesize#1};
		\node[below,inline text] at (0,-.5) {\footnotesize#2};
		\node[above,inline text]  at (0,.5) {\footnotesize#3};
}
\def\Raction[#1][#2][#3]{
	\draw (0,-.5)--(0,.5);
	\draw (.6,-.5)to[out=90,in=-30](0,0);
	\node[below,inline text] at (.6,-.5) {\footnotesize#1};
	\node[below,inline text] at (0,-.5) {\footnotesize#2};
	\node[above,inline text] at (0,.5) {\footnotesize#3};
}

\def\Lassociator[#1][#2][#3][#4]{
	\draw (0,-.5)--(0,.5);
	\draw (-.6,-.5)to[out=90,in=210](0,0);
	\draw (.6,-.5)to[out=90,in=-30](0,.25);
	\node[below,inline text] at (-.6,-.5) {\footnotesize#1};
	\node[below,inline text] at (.6,-.5) {\footnotesize#2};
	\node[below,inline text] at (0,-.5) {\footnotesize#3};
	\node[above,inline text] at (0,.5) {\footnotesize#4};
}
\def\Rassociator[#1][#2][#3][#4]{
	\draw (0,-.5)--(0,.5);
	\draw (-.6,-.5)to[out=90,in=210](0,.25);
	\draw (.6,-.5)to[out=90,in=-30](0,0);
	\node[below,inline text] at (-.6,-.5) {\footnotesize#1};
	\node[below,inline text] at (.6,-.5) {\footnotesize#2};
	\node[below,inline text] at (0,-.5) {\footnotesize#3};
	\node[above,inline text] at (0,.5) {\footnotesize#4};
}

\usepackage{hhline}
\usepackage{etoolbox}

\let\originalleft\left
\let\originalright\right
\renewcommand{\left}{\mathopen{}\mathclose\bgroup\originalleft}
\renewcommand{\right}{\aftergroup\egroup\originalright}

\newcommand{\vvec}[1]
{
	\operatorname{\bf Vec}
	\ifstrequal{#1}{}
	{}
	{\left(#1\right)}
}

\newcommand{\vvectwist}[2]
{
	\operatorname{\bf Vec}^{#2}
	\ifstrequal{#1}{}
	{}
	{\left(#1\right)}
}

\newcommand{\lad}[2]
{
	\operatorname{\bf Lad}
	\ifstrequal{#2}{}
	{}
	{_{#2}}
	\ifstrequal{#1}{}
	{}
	{\left(#1\right)}
}

\newcommand{\kar}[1]
{
	\operatorname{\bf Kar}
	\ifstrequal{#1}{}
	{}
	{\left(#1\right)}
}

\newcommand{\bpr}[1]
{
	\operatorname{\bf BPR}
	\ifstrequal{#1}{}
	{}
	{\left(#1\right)}
}

\newcommand{\dih}[1]
{
	\operatorname{Dih}
	\ifstrequal{#1}{}
	{}
	{_{#1}}
}

\newcommand{\trivalentvertex}[4]{
   \draw (0,0)--(0,1);
   \node[above,inline text] at (0,1) {#1};
   \draw (0,0)--(-0.707,-0.707);
   \node[below,inline text] at (-0.707,-0.707) {#2};
   \draw (0,0)--(0.707,-0.707);
   \node[below,inline text] at (0.707,-0.707) {#3};
   \node[left,inline text] at (0,0) {#4};
}

\usepackage{soul}

\allowdisplaybreaks
\begin{document}

\title{Domain walls in topological phases and the Brauer-Picard ring for $\vvec{\ZZ{p}}$}

\author{Daniel Barter}
\email{danielbarter@gmail.com}
\affiliation{Mathematical Sciences Institute, Australian National University, Canberra, Australia}
\author{Jacob C.\ Bridgeman}
\email{jcbridgeman1@gmail.com}
\affiliation{Centre for Engineered Quantum Systems, School of Physics, The University of Sydney, Sydney, Australia}
\author{Corey Jones}
\email{cormjones88@gmail.com} 
\affiliation{Mathematical Sciences Institute, Australian National University, Canberra, Australia}

\date{\today}

\begin{abstract}
	We show how to calculate the relative tensor product of bimodule categories (not necessarily invertible) using ladder string diagrams. As an illustrative example, we compute the Brauer-Picard ring for the fusion category $\vvec{\ZZ{p}}$. Moreover, we provide a physical interpretation of all indecomposable bimodule categories in terms of domain walls in the associated topological phase. We show how this interpretation can be used to compute the Brauer-Picard ring from a physical perspective.
\end{abstract}

\maketitle

%
\vspace*{3mm}
\section{Introduction}\label{S:introduction}
One of the most important challenges on the path to building a quantum computer is choosing a substrate for encoding the quantum information. Such a substrate should allow for robust storage and manipulation of the qubits, even in the presence of environmental noise. These requirements have prompted research into topological phases of matter\cite{Dennis2002,MR1951039,Brown2014,Terhal2015}. These quantum systems are formed from many local degrees of freedom which interact in such a way as to produce emergent, global, topological properties. In this paper, we consider domain walls between 2-dimensional, non-chiral topological phases with long-range entanglement. 

In addition to the bulk properties of topological phases, it has become clear that quantum information can be encoded and manipulated using domain walls, and defects thereof, between bulks\cite{Dennis2002,0610082,Bombin2007a,Bombin2010,Brown2013a,Pastawski2014,Yoshida2015a,Brown2016,IrisCong1,IrisCong2,Yoshida2017,Brown2018}. The inclusion of domain walls, both invertible and non-invertible, can increase both the encoding rate and the number of fault-tolerant quantum operations that can be performed. It is therefore important to study the full set of possible domain walls for a given topological phase.

The renormalization invariant properties of such topological phases are described by a (2+1)D topological quantum field theory (TQFT). Depending on the application, there are several different descriptions of these TQFTs. The field-theoretic approach due to Witten\cite{MR990772} uses a Lagrangian which does not depend on the space-time metric. In mathematics, (2+1)D TQFTs are described using functors from a cobordism category into some algebraic category, and can be used to obtain invariants of framed 3-manifolds\cite{MR1001453,MR1355899}. A vast amount of work has been done explaining how these functors can be constructed directly using fusion categories. We refer the reader to \onlinecite{MR3674995} for a good survey of the literature. Most relevant to our current motivation, in condensed matter physics and quantum information theory, one often proceeds by constructing an explicit microscopic model which realizes a given TQFT in the low energy space. One of the most famous and important models of this kind is known as the toric- or surface- code\cite{MR1951039}. In addition to its theoretical importance, this model is of great experimental interest\cite{chow2014implementing,Gambetta1}. More generally, given a fusion category, the Levin-Wen model\cite{Levin2005} is a lattice Hamiltonian that realizes the associated topological phase.

\begin{figure}[t]
	$\begin{array}{c}
	\includeTikz{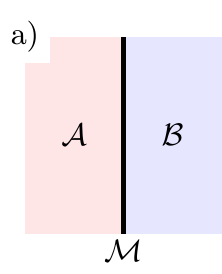}
	{
		\begin{tikzpicture}
		\fill[red!10](-1,-1) rectangle (0,1);
		\fill[blue!10](1,-1) rectangle (0,1);
		\draw[ultra thick] (0,-1)--(0,1);
		\node at (-.5,0) {$\mathcal{A}$};
		\node at (.5,0) {$\mathcal{B}$};
		\node[below,inline text] at (0,-1) {$\mathcal{M}$};
		\node[above,inline text,white] at (0,1) {$\mathcal{M}$};
		\node[fill=white] at (-1,1) {a)};
		\end{tikzpicture}
	}
	\end{array}$
	\hspace{15mm}
	$\begin{array}{c}
	\includeTikz{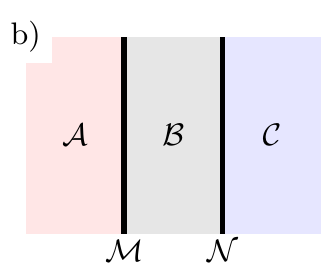}
	{
		\begin{tikzpicture}
		\fill[red!10](-1,-1) rectangle (0,1);
		\fill[black!10](1,-1) rectangle (0,1);
		\fill[blue!10](1,-1) rectangle (2,1);
		\draw[ultra thick] (0,-1)--(0,1);
		\draw[ultra thick] (1,-1)--(1,1);
		\node at (-.5,0) {$\mathcal{A}$};
		\node at (.5,0) {$\mathcal{B}$};
		\node at (1.5,0) {$\mathcal{C}$};
		\node[below,inline text] at (0,-1) {$\mathcal{M}$};
		\node[below,inline text] at (1,-1) {$\mathcal{N}$};
		\node[above,inline text,white] at (0,1) {$\mathcal{M}$};
		\node[above,inline text,white] at (1,1) {$\mathcal{N}$};
		\node[fill=white] at (-1,1) {b)};
		\end{tikzpicture}
	}
	\end{array}
	=\sum
	\begin{array}{c}
	\includeTikz{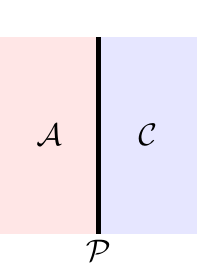}
	{
		\begin{tikzpicture}
		\fill[red!10](-1,-1) rectangle (0,1);
		\fill[blue!10](1,-1) rectangle (0,1);
		\draw[ultra thick] (0,-1)--(0,1);
		\node at (-.5,0) {$\mathcal{A}$};
		\node at (.5,0) {$\mathcal{C}$};
		\node[below,inline text] at (0,-1) {$\mathcal{P}$};
		\node[above,inline text,white] at (0,1) {$\mathcal{P}$};
		\end{tikzpicture}
	}
	\end{array}$
	\caption{a) Domain walls between topological phases corresponding to fusion categories $\mathcal{A}$ and $\mathcal{B}$ are described by a bimodule category $\mathcal{A}\curvearrowright\mathcal{M}\curvearrowleft\mathcal{B}$.\\
		\hspace*{11mm} b) In this paper we discuss the concatenation of domain walls, corresponding to the relative tensor product of bimodules $\mathcal{M}\otimes_\mathcal{B}\mathcal{N}$.} \label{fig:combiningDW}
\end{figure}

In terms of fusion categories, domain walls with no terminating defects are described by module categories. More precisely, if $\mathcal{A}$ and $\mathcal{B}$ are fusion categories, then a domain wall between the corresponding topological phases can be described using a bimodule category $\mathcal{A} \curvearrowright \mathcal{M} \curvearrowleft \mathcal{B}$\cite{FUCHS2002353,MR3370609,MR2942952,Kong2013} as shown in Fig.~\hyperref[fig:combiningDW]{\ref*{fig:combiningDW}a}. The bimodule encodes the information required to glue the two bulks together. It is natural to consider combining adjacent domain walls, and considering them as a single interface (Fig.~\hyperref[fig:combiningDW]{\ref*{fig:combiningDW}b}). Mathematically, this corresponds to the relative tensor product of bimodules, $\mathcal{M} \otimes_{\mathcal{B}} \mathcal{N}$.

The collection of $\mathcal{C}$-\hspace{0.02cm}$\mathcal{C}$ bimodules equipped with the relative tensor product $\otimes_{\mathcal{C}}$ is called the Brauer-Picard ring of $\mathcal{C}$, which we denote $\bpr{\mathcal{C}}$. It has a basis given by equivalence classes of the indecomposable bimodules. The group of units of $\bpr{\mathcal{C}}$ is known as the Brauer-Picard group. These transparent domain walls are closely related to the locality preserving unitary transformations that can be applied to the associated quantum code in a fault-tolerant manner\cite{Yoshida2017}. Inclusion of non-invertible elements allows for a richer theory including punctures and more general objects\cite{Brown2016}.

The mathematical question at the heart of this paper is:
\begin{question}
	Compute the Brauer-Picard ring for a semi-simple fusion category.
\end{question}
Although we cannot provide a complete answer to this question, we present a general computational method. We apply this method to answer the question:
\begin{question}\label{questionB}
	Compute the Brauer-Picard ring for $\vvec{\ZZ{p}}$ with $p$ prime.
\end{question}

In this paper, we use \emph{ladder string diagrams}\cite{MR2978449,MR2187404} to answer Question~\ref{questionB}.
These ladder string diagrams are not new. They were first introduced by Morrison and Walker in Chapter~6 of \onlinecite{MR2978449} as a tool for computing the relative tensor products for modules over $n$-categories. The string diagram calculus for fusion categories as described in \onlinecite{MR2767048}, allows us to interpret fusion categories as disc-like 2-categories, and use the ideas of Morrison and Walker to compute relative tensor products.

The problem of computing the Brauer-Picard ring of a fusion category is not new. In \onlinecite{MR2677836}, a formula for the Brauer-Picard ring of $\vvec{A}$ when $A$ is an Abelian group is given. Their formula gives the same table as ours, which was verified in \onlinecite{MR3370609} for the invertible bimodules and in \onlinecite{1809.00245} for $\vvec{\ZZ{2}}$. Our computation method is completely independent of the formula in \onlinecite{MR2677836}, and generalizes to nonabelian groups. This will be the topic of future work\footnote{D.~Barter, J.~C.~Bridgeman and C.~Jones, \emph{in preparation}}.

We give physical interpretations for all the indecomposable bimodules and the relative tensor product multiplication table. In addition to providing an illustrative example of the more general computation, this class of models includes the toric code (associated to $\vvec{\ZZ{2}}$). This model is the most studied topological phase, both as a condensed matter model, and in relation to quantum information processing.

This paper is structured as follows: In Section~\ref{S:prelim}, we provide some mathematical preliminaries that are required for the remainder of the paper. We then describe the general computation of the relative tensor product of bimodules. We also prove that the category of representations of the ladder category corresponds to the relative tensor product. In Section~\ref{S:Zp}, we first explain the well known classification of indecomposable bimodules for the fusion category $\vvec{\ZZ{p}}$. We then provide several detailed computations of the relative tensor product using ladder string diagrams. The section concludes with a physical interpretation of the bimodules and their multiplication table. In Section~\ref{S:conc} we summarize and discuss some extensions of this work.

%

\section{Preliminaries} \label{S:prelim}

Throughout this paper, all categories are enriched over $\vvec{}$, the category of $\mathbb{C}$-vector spaces. This means that the Hom spaces are vector spaces and composition is bilinear. In particular, the endomorphisms of any object form an algebra. We repeat some standard notation for tensor categories, as found in \onlinecite{MR3242743}.

\begin{definition} \label{deligne_tensor_product}
  If $\mathcal{C}$ and $\mathcal{D}$ are semi-simple categories enriched over ${\bf Vec}$, then we can form the {\em Deligne tensor product} $\mathcal{C} \otimes \mathcal{D}$ as defined on Page~15 of \onlinecite{MR3242743}. If $c \in \mathcal{C}$ and $d \in \mathcal{D}$ are objects, then $(c,d)$ is an object in $\mathcal{C} \otimes \mathcal{D}$. The Hom spaces are defined by
  \begin{align}(\mathcal{C} \otimes \mathcal{D})\left((c_1,d_1),(c_2,d_2)\right) = \mathcal{C}(c_1,c_2) \otimes \mathcal{D}(d_1,d_2),\end{align}
  where the right hand side is the usual tensor product of vector spaces.
  The Deligne tensor product is the correct notion of product for categories enriched over vector spaces.
\end{definition}

\begin{definition} \label{tensor_cat_def} A {\em tensor category} $\mathcal{C}$ is a category $\mathcal{C}$ equipped with a functor $- \otimes - : \mathcal{C} \otimes \mathcal{C} \to \mathcal{C}$, a natural isomorphism $ (- \otimes -) \otimes - \cong - \otimes (- \otimes -)$ called the associator and a special object $1 \in \mathcal{C}$ which satisfy the pentagon equation and unit equations respectively. These can be found on Page~22 of \onlinecite{MR3242743}. If $\mathcal{C}$ is semi-simple, using the string diagram notation as explained in \onlinecite{MR2767048}, a vector in $\mathcal{C}(a \otimes b, c)$ can be represented by a trivalent vertex: 
\begin{align}
\begin{array}{c}
  \includeTikz{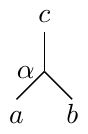}{
    \begin{tikzpicture}[scale=0.4, every node/.style={scale=0.8}]
      \trivalentvertex{$c$}{$a$}{$b$}{$\alpha$};
    \end{tikzpicture}
  }
\end{array}.
\end{align}

If we choose bases for all the vector spaces $\mathcal{C}(a \otimes b,c)$, then the associator can be represented as a tensor F, where
  \begin{align}
    \begin{array}{c}
	    \includeTikz{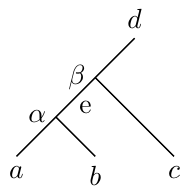}
	    {
	    	\begin{tikzpicture}[scale=0.4, every node/.style={scale=0.8}]
		    	\draw(0,0)--(3,3) node[below] at (0,0) {$a$} node[above] at (3,3) {$d$};
		    	\draw(2,0)--(1,1) node[below] at (2,0) {$b$};
		    	\draw(4,0)--(2,2) node[below] at (4,0) {$c$};
		    	\node[] at (1.75,1.25) {e};
		    	\node[left] at(1,1){$\alpha$};
		    	\node[left] at(2,2){$\beta$};
	    	\end{tikzpicture}
	    }
    \end{array}
    =\sum_{(f,\mu,\nu)}
    \biggl[F^{abc}_{d}\biggr]_{(e,\alpha,\beta)(f,\mu,\nu)}
    \begin{array}{c}
    \includeTikz{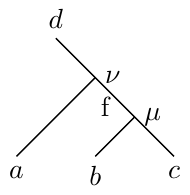}
    {
    	\begin{tikzpicture}[scale=0.4, every node/.style={scale=0.8}]
    	\draw(0,0)--(-3,3) node[below] at (0,0) {$c$} node[above] at (-3,3) {$d$};
    	\draw(-2,0)--(-1,1) node[below] at (-2,0) {$b$};
    	\draw(-4,0)--(-2,2) node[below] at (-4,0) {$a$};
    	\node at (-1.75,1.25) {f};
    	\node[right] at(-1,1){$\mu$};
    	\node[right] at(-2,2){$\nu$};
    	\end{tikzpicture}
    }
    \end{array}.
  \end{align}
  where $\alpha \in \mathcal{C}(a \otimes b,e), \beta \in \mathcal{C}(e \otimes c,d), \mu \in \mathcal{C}(b \otimes c,f)$ and $\nu \in \mathcal{C}(a \otimes f,d)$ are basis vectors.
\end{definition}

\begin{definition} \label{fusion_cat_def} A {\em fusion category} is a semi-simple rigid tensor category $\mathcal{C}$ with a finite number of simple objects and a simple unit. An object is {\em simple} if there are no non-trivial sub-objects. Semi-simple means that every object in $\mathcal{C}$ is a direct sum of simple objects. Rigid is a technical condition defined on Page~40 of \onlinecite{MR3242743} which implies that objects in $\mathcal{C}$ have duals and they behave like duals in the category of vector spaces. When doing fusion category computations, it is customary to fix a set of simple objects and express all other objects as direct sums of the chosen simple ones. 
\end{definition}

\begin{definition} \label{module_cat_def}
  If $\mathcal{C}$ is a fusion category, then a {\em module category} $\mathcal{M}$ is a semi-simple category equipped with a functor $ - \vartriangleright - : \mathcal{C} \otimes \mathcal{M} \to \mathcal{M}$ and a natural isomorphism $ - \vartriangleright ( - \vartriangleright - ) \cong (- \otimes - ) \vartriangleright -$ which satisfies some coherence conditions that can be found on Page~131 of \onlinecite{MR3242743}. Using the string diagram notation, if we choose basis vectors for all the vector spaces $\mathcal{M}(a \vartriangleright m, n)$ then the associator can be represented using a tensor
  \begin{align}
  \begin{array}{c}
	  \includeTikz{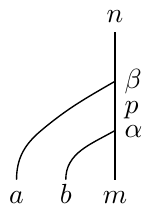}
	  {
	  	\begin{tikzpicture}[scale=0.5, every node/.style={scale=0.8}]
	  	\draw(0,0)--(0,3);
	  	\draw(-1,0)to[out=90,in=210](0,1);
	  	\draw(-2,0)to[out=90,in=220](-1.5,1)to[out=40,in=210](0,2);
	  	\node[below,inline text] at (-2,0) {$a$};
	  	\node[below,inline text] at (-1,0) {$b$};
	  	\node[below,inline text] at (0,0) {$m$};
	  	\node[above,inline text] at (0,3) {$n$};
	  	\node[right,inline text] at (0,1) {$\alpha$};
	  	\node[right,inline text] at (0,1.5) {$p$};
	  	\node[right,inline text] at (0,2) {$\beta$};
	  	\end{tikzpicture}
	  }
  \end{array}
  =\sum_{(q,\mu,\nu)}
  \biggl[L^{abm}_{n}\biggr]_{(p,\alpha,\beta)(q,\mu,\nu)}
  \begin{array}{c}
  \includeTikz{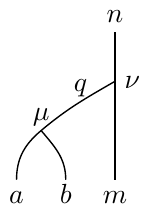}
  {
  	\begin{tikzpicture}[scale=0.5, every node/.style={scale=0.8}]
	  	\draw(0,0)--(0,3);
	  	\draw(-1,0)to[out=90,in=-50](-1.5,1);
	  	\draw(-2,0)to[out=90,in=220](-1.5,1)to[out=40,in=210](0,2);
	  	\node[below,inline text] at (-2,0) {$a$};
	  	\node[below,inline text] at (-1,0) {$b$};
	  	\node[below,inline text] at (0,0) {$m$};
	  	\node[above,inline text] at (0,3) {$n$};
	  	\node[above,inline text] at (-1.5,1) {$\mu$};
	  	\node[above,inline text] at (-.7,1.6) {$q$};
	  	\node[right,inline text] at (0,2) {$\nu$};
  	\end{tikzpicture}
  }
  \end{array}.\label{eqn:leftassociatorA}
  \end{align}
  We shall often abuse notation and use the $\otimes$ symbol instead of $\vartriangleright$. The functor can always be inferred by scrutinizing the category in which the left object lives. 
\end{definition}

\begin{definition} \label{bimodule_definition}
  Let $\mathcal{C},\,\mathcal{D}$ be fusion categories. A {\em bimodule category} $\mathcal{C} \curvearrowright \mathcal{M} \curvearrowleft \mathcal{D}$ is a semi-simple category equipped with functors $ - \vartriangleright - : \mathcal{C} \otimes \mathcal{M} \to \mathcal{M}$ and $ - \vartriangleleft - : \mathcal{M} \otimes \mathcal{D} \to \mathcal{M}$ and three natural isomorphisms
  \begin{align}
    &- \vartriangleright (- \vartriangleright -) \cong (- \otimes - ) \vartriangleright - \\
    & (- \vartriangleleft -) \vartriangleleft - \cong - \vartriangleleft (- \otimes -) \\
    & - \vartriangleright (- \vartriangleleft -) \cong (- \vartriangleright -) \vartriangleleft -.
  \end{align}
  If we choose bases for the Hom spaces $\mathcal{M}(c \vartriangleright m, n)$ and $\mathcal{M}(m \vartriangleleft d,n)$, then these natural isomorphisms can be represented as tensors using Equation~\ref{eqn:leftassociatorA} and
  
  \begin{align}
  \begin{array}{c}
  \includeTikz{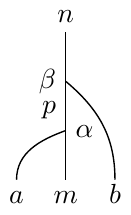}
  {
  	\begin{tikzpicture}[scale=0.5, every node/.style={scale=0.8}]
  	\draw(0,0)--(0,3);
  	\draw(1,0)to[out=90,in=-40](0,2);
  	\draw(-1,0)to[out=90,in=200](0,1);
  	\node[below,inline text] at (-1,0) {$a$};
  	\node[below,inline text] at (1,0) {$b$};
  	\node[below,inline text] at (0,0) {$m$};
  	\node[above,inline text] at (0,3) {$n$};
  	\node[right,inline text] at (0,1) {$\alpha$};
  	\node[left,inline text] at (0,1.5) {$p$};
  	\node[left,inline text] at (0,2) {$\beta$};
  	\end{tikzpicture}
  }
  \end{array}
  &=\sum_{(q,\mu,\nu)}
  \biggl[C^{abm}_{n}\biggr]_{(p,\alpha,\beta)(q,\mu,\nu)}
  \begin{array}{c}
  \includeTikz{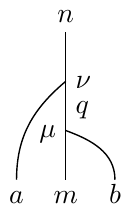}
  {
  	\begin{tikzpicture}[scale=0.5, every node/.style={scale=0.8}]
  	\draw(0,0)--(0,3);
  	\draw(1,0)to[out=90,in=-20](0,1);
  	\draw(-1,0)to[out=90,in=220](0,2);
  	\node[below,inline text] at (-1,0) {$a$};
  	\node[below,inline text] at (1,0) {$b$};
  	\node[below,inline text] at (0,0) {$m$};
  	\node[above,inline text] at (0,3) {$n$};
  	\node[left,inline text] at (0,1) {$\mu$};
  	\node[right,inline text] at (0,1.5) {$q$};
  	\node[right,inline text] at (0,2) {$\nu$};
  	\end{tikzpicture}
  }
  \end{array},\\
  \begin{array}{c}
  \includeTikz{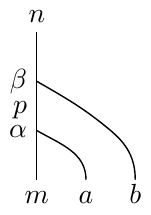}
  {
  	\begin{tikzpicture}[xscale=-0.5,yscale=.5, every node/.style={scale=0.8}]
  	\draw(0,0)--(0,3);
  	\draw(-1,0)to[out=90,in=210](0,1);
  	\draw(-2,0)to[out=90,in=220](-1.5,1)to[out=40,in=210](0,2);
  	\node[below,inline text] at (-2,0) {$b$};
  	\node[below,inline text] at (-1,0) {$a$};
  	\node[below,inline text] at (0,0) {$m$};
  	\node[above,inline text] at (0,3) {$n$};
  	\node[left,inline text] at (0,1) {$\alpha$};
  	\node[left,inline text] at (0,1.5) {$p$};
  	\node[left,inline text] at (0,2) {$\beta$};
  	\end{tikzpicture}
  }
  \end{array}
  &=\sum_{(q,\mu,\nu)}
  \biggl[R^{mab}_{n}\biggr]_{(p,\alpha,\beta)(q,\mu,\nu)}
  \begin{array}{c}
  \includeTikz{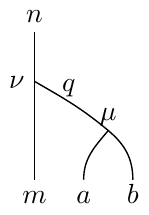}
  {
  	\begin{tikzpicture}[xscale=-0.5,yscale=.5, every node/.style={scale=0.8}]
  	\draw(0,0)--(0,3);
  	\draw(-1,0)to[out=90,in=-50](-1.5,1);
  	\draw(-2,0)to[out=90,in=220](-1.5,1)to[out=40,in=210](0,2);
  	\node[below,inline text] at (-2,0) {$b$};
  	\node[below,inline text] at (-1,0) {$a$};
  	\node[below,inline text] at (0,0) {$m$};
  	\node[above,inline text] at (0,3) {$n$};
  	\node[above,inline text] at (-1.5,1) {$\mu$};
  	\node[above,inline text] at (-.7,1.6) {$q$};
  	\node[left,inline text] at (0,2) {$\nu$};
  	\end{tikzpicture}
  }
  \end{array}.
  \end{align}
There is an equivalent way to present bimodules. Let $\mathcal{D}^{\rm mp}$ be the fusion category which is $\mathcal{D}$ as a category, but with the reversed tensor product\footnote{There are multiple ways to take the opposite of a tensor category. The reader should consult \onlinecite{MR3242743} for the definitions of all the tensor category opposite constructions and how they are related}. Then a bimodule $\mathcal{C} \curvearrowright \mathcal{M} \curvearrowleft \mathcal{D}$ is the same data as a module $\mathcal{C} \otimes \mathcal{D}^{\rm mp} \curvearrowright \mathcal{M}$. As in Definition~\ref{module_cat_def}, we shall often abuse notation and write $\otimes$ instead of $\vartriangleright$ and $\vartriangleleft$.
\end{definition}

\begin{definition} \label{relative_tensor_products}
   Let $\mathcal{B}$ be a fusion category and suppose we have modules $\mathcal{M} \curvearrowleft \mathcal{B} \curvearrowright \mathcal{N}$. Then the relative tensor product $\mathcal{M} \otimes_{\mathcal{B}} \mathcal{N}$ is defined as follows: the objects consist of objects $(m,n) \in \mathcal{M} \otimes \mathcal{N}$ together with specified isomorphisms $(m \otimes b,n) \cong (m,b \otimes n)$ for all $b \in \mathcal{B}$ that satisfy some coherence conditions. The morphisms in $\mathcal{M} \otimes_{\mathcal{B}} \mathcal{N}$ are morphisms in $\mathcal{M} \otimes \mathcal{N}$ that are compatible with this data. This definition appears in \onlinecite{1406.4204}. 
\end{definition}

\begin{definition} \label{ladder_cat_def}
  Let $\mathcal{B}$ be a fusion category and suppose we have modules $\mathcal{M} \curvearrowleft \mathcal{B} \curvearrowright \mathcal{N}$. The {\em ladder category} $\lad{\mathcal{M},\mathcal{N}}{\mathcal{B}}$ is defined as follows: The objects are pairs of simple objects $(m \in \mathcal{M},n \in \mathcal{N})$ and the morphisms are linear combinations of basic ladder string diagrams
\begin{align}
\begin{array}{c}
\includeTikz{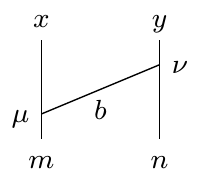}
{
	\begin{tikzpicture}
		\coordinate (A) at (0,0);
		\ladder[A][$m$][$n$][$b$][$x$][$y$];
		\node[inline text,left] at (-.6,-.25) {\footnotesize$\mu$};
		\node[inline text,right] at (.6,.25) {\footnotesize$\nu$};
	\end{tikzpicture}
}
\end{array},
\end{align}
built out of trivalent vertices in $\mathcal{M}(m,x \otimes b)$ and $\mathcal{N}(b \otimes n,y)$. Composition is given by stacking ladders. We can use the decomposition
\begin{align}
\begin{array}{c}
\includeTikz{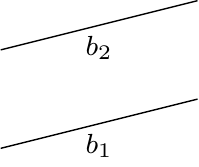}
	{
	\begin{tikzpicture}
		\draw(-1,-.75)--(1,-.25) node[midway,inline text,below] {\footnotesize$b_1$};
		\draw(-1,.25)--(1,.75) node[midway,inline text,below] {\footnotesize$b_2$};
	\end{tikzpicture}
	}
\end{array}
&=
\sum_{b,\alpha}
\begin{array}{c}
\includeTikz{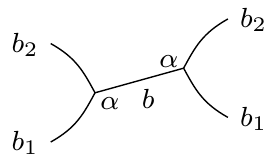}
{
	\begin{tikzpicture}[xscale=1.5]
		\coordinate (A) at (0,0);
		\coordinate (B) at (0,1);
		\draw (-.6,0)to[out=40,in=250](-.3,.5);
		\draw[yscale=-1] (-.6,-1)to[out=40,in=250](-.3,-.5);
		\draw(-.3,.5)--(.3,.75) node[pos=.6,inline text,below] {\footnotesize$b$};
		\begin{scope}[xscale=-1,shift={(0,.25)}]
		\draw (-.6,0)to[out=40,in=250](-.3,.5);
		\draw[yscale=-1] (-.6,-1)to[out=40,in=250](-.3,-.5);
		\node[inline text] at (-.2,.6) {\footnotesize$\alpha$};
		\end{scope}
		\node[inline text] at (-.2,.42) {\footnotesize$\alpha$};
		\node[inline text,left] at (-.6,0) {\footnotesize$b_1$};
		\node[inline text,right] at (.6,.25) {\footnotesize$b_1$};
		\node[inline text,left] at (-.6,1) {\footnotesize$b_2$};
		\node[inline text,right] at (.6,1.25) {\footnotesize$b_2$};
	\end{tikzpicture}
}
\end{array},
\end{align}
to compute the resulting morphism
\begin{align}
\begin{array}{c}
\includeTikz{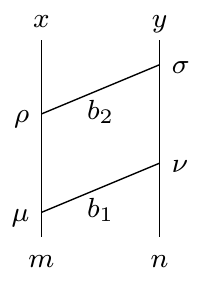}
{
	\begin{tikzpicture}
	\coordinate (A) at (0,0);
	\coordinate (B) at (0,1);
	\ladder[A][$m$][$n$][$b_1$][][];
	\ladder[B][][][$b_2$][$x$][$y$];
		\node[inline text,left] at (-.6,-.25) {\footnotesize$\mu$};
		\node[inline text,right] at (.6,.25) {\footnotesize$\nu$};
		\node[inline text,left] at (-.6,.75) {\footnotesize$\rho$};
		\node[inline text,right] at (.6,1.25) {\footnotesize$\sigma$};
	\end{tikzpicture}
}
\end{array}&=\sum_{b,\alpha}
\begin{array}{c}
	\includeTikz{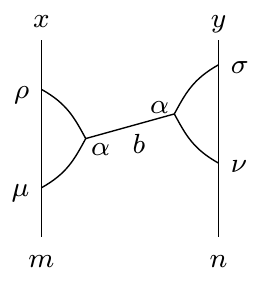}
	{
		\begin{tikzpicture}[xscale=1.5]
		\coordinate (A) at (0,0);
		\coordinate (B) at (0,1);
		\ladder[A][$m$][$n$][][][];
		\ladder[B][][][][$x$][$y$];
		\draw (-.6,0)to[out=40,in=250](-.3,.5);
		\draw[yscale=-1] (-.6,-1)to[out=40,in=250](-.3,-.5);
		\draw(-.3,.5)--(.3,.75) node[pos=.6,inline text,below] {\footnotesize$b$};
		\begin{scope}[xscale=-1,shift={(0,.25)}]
		\draw (-.6,0)to[out=40,in=250](-.3,.5);
		\draw[yscale=-1] (-.6,-1)to[out=40,in=250](-.3,-.5);
		\node[inline text] at (-.2,.6) {\footnotesize$\alpha$};
		\end{scope}
		\node[inline text] at (-.2,.42) {\footnotesize$\alpha$};
		\node[inline text,left] at (-.6,0) {\footnotesize$\mu$};
		\node[inline text,right] at (.6,.25) {\footnotesize$\nu$};
		\node[inline text,left] at (-.6,1) {\footnotesize$\rho$};
		\node[inline text,right] at (.6,1.25) {\footnotesize$\sigma$};
		\end{tikzpicture}
	}
\end{array}.\label{eqn:stackingladders}
\end{align}
In the module category, the triangles reduce to a sum of basis trivalent vertices
\begin{align}
\begin{array}{c}
\includeTikz{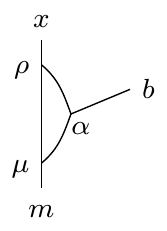}{
	\begin{tikzpicture}
	\draw (-.6,-.25)--(-.6,1.25) node[pos=0,inline text,below] {\footnotesize$m$} node[pos=1,inline text,above] {\footnotesize$x$};
	\draw (-.6,0)to[out=40,in=250](-.3,.5);
	\draw[yscale=-1] (-.6,-1)to[out=40,in=250](-.3,-.5);
	\draw(-.3,.5)--(.3,.75) node[pos=1,inline text,right] {\footnotesize$b$};
	\node[inline text,below] at (-.2,.6) {\footnotesize$\alpha$};
		\node[inline text,left] at (-.6,0) {\footnotesize$\mu$};
		\node[inline text,left] at (-.6,1) {\footnotesize$\rho$};
	\end{tikzpicture}}
\end{array}
&=
\sum_{\beta}
c_\beta
\begin{array}{c}
\includeTikz{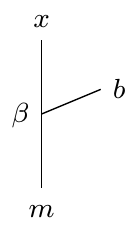}{
	\begin{tikzpicture}
	\draw (-.3,-.25)--(-.3,1.25) node[pos=0,inline text,below] {\footnotesize$m$} node[pos=1,inline text,above] {\footnotesize$x$};
	\draw(-.3,.5)--(.3,.75) node[pos=1,inline text,right] {\footnotesize$b$};
	\node[inline text,left] at (-.3,.5) {\footnotesize$\beta$};
	\end{tikzpicture}}
\end{array},
\end{align}
which can be used to reduce Eqn.~\ref{eqn:stackingladders} to a linear combination of ladders. The coefficient $c_\beta$ is determined by the choice of basis trivalent vertices.

\end{definition}

\begin{definition}[Karoubi envelope] \label{Karoubi_envelope_definition}
Let $\mathcal{C}$ be a category. As explained in \onlinecite{MR2253455}, the {\em Karoubi envelope} $\kar{\mathcal{C}}$ has objects $(A,e)$ where $A \in \mathcal{C}$ is an object and $e : A \to A$ is an idempotent. A morphism $f : (A,e) \to (A',e')$ is a morphism $f : A \to A'$ in $\mathcal{C}$ such that $e' f = f e = f$. If $\mathcal{C}$ is semi-simple, it follows from Yoneda's lemma that $\kar{\mathcal{C}}$ is equivalent to the category $[\mathcal{\mathcal{C}}^{\rm op},\vvec{}]$, whose objects are functors $\mathcal{C}^{\rm op} \to {\bf Vec}$ and morphisms the natural transformations. Under this equivalence, the idempotent $e : A \to A$ is sent to the functor $e\mathcal{C}(-,A)$: morphisms into $A$ which factor through $e$.
\end{definition}
\begin{definition} \label{vec_G_def} 
  Let $G$ be a finite group. If $V$ is a vector space, a $G$-grading on $V$ is a direct sum decomposition $V = \oplus_{g \in G} V_g$. We define $\vvec{G}$ to be the category of $G$-graded vector spaces and linear maps which preserve the grading. The tensor product of $G$-graded vector spaces is defined by
  \begin{align}(V \otimes W)_k = \bigoplus_{gh = k} V_g \otimes W_h.\end{align}
The category $\vvec{G}$ is semi-simple and the simple objects are parameterized by $G$. The element $g \in G$ corresponds to the $G$-graded vector space which is 1 dimensional in degree $g$ and 0 dimensional in all other degrees. We shall identify the simple object corresponding to $g \in G$ with the element $g$ itself. Using this notation, the tensor product is exactly $g \otimes h = gh$.
\end{definition}

Let $\mathcal{A},\,\mathcal{B}$ and $\mathcal{C}$ be fusion categories. Suppose that we have bimodules $\mathcal{A} \curvearrowright \mathcal{M} \curvearrowleft \mathcal{B} \curvearrowright \mathcal{N} \curvearrowleft \mathcal{C}$. Our main goal in this paper is to compute the relative tensor product $\mathcal{M} \otimes_{\mathcal{B}} \mathcal{N}$. The category $\mathcal{M} \otimes_{\mathcal{B}} \mathcal{N}$ inherits an $\mathcal{A}$-$\mathcal{C}$ bimodule structure, so can be interpreted as a domain wall. By computing the the relative tensor product, we can decompose the concatenated domain wall into a sum of indecomposable domain walls: ${\rm wall}_\mathcal{M}\times{\rm wall}_\mathcal{N}=\sum_\mathcal{P}{\rm wall}_\mathcal{P}$. Our main computational tool is the following:
\begin{proposition} \label{ladder_model_relative_tensor_product}
  Let $\mathcal{B}$ be a fusion category and suppose that we have module categories $\mathcal{M} \curvearrowleft \mathcal{B} \curvearrowright \mathcal{N}$. Then $ \mathcal{M} \otimes_{\mathcal{B}} \mathcal{N}$ is equivalent to representations of $\lad{\mathcal{M},\mathcal{N}}{\mathcal{B}}^{\rm op}$. If $\mathcal{C}$ is semi-simple, we get the equivalence
\begin{align}\mathcal{M} \otimes_{\mathcal{B}} \mathcal{N} \cong \kar{\protect\lad{\mathcal{M},\mathcal{N}}{\mathcal{B}}}.\end{align}
\end{proposition}

\noindent Proposition~\ref{ladder_model_relative_tensor_product} is the result which lets us compute relative tensor products. It is well known to mathematicians who work with fusion categories. It first appeared in Chapter~6 of \onlinecite{MR2978449} as a tool for computing the relative tensor products for modules over $n$-categories. We include a proof of Proposition~\ref{ladder_model_relative_tensor_product} specialized to our setting. Our proof uses the following proposition:
\begin{proposition} \label{eno_remark_3.9}
Let $\mathcal{A},\mathcal{B},\mathcal{C}$ be fusion categories and $\mathcal{A} \curvearrowright \mathcal{M} \curvearrowleft \mathcal{B} \curvearrowright \mathcal{N} \curvearrowleft \mathcal{C}$ be bimodules. Let $\Lambda$ be the set of isomorphism classes of simple objects in $\mathcal{B}$. Then
\begin{align}
I &= \bigoplus_{b \in \Lambda} b^* \otimes b
\end{align}
has a canonical algebra structure $\mathcal{B}^{\otimes {\rm op}} \otimes \mathcal{B}$ and $\mathcal{M} \otimes_{\mathcal{B}} \mathcal{N}$ is equivalent to the category ${\bf Mod}(I)$ of $I$-modules in $\mathcal{M} \otimes \mathcal{N}$. 
\end{proposition}
\noindent This is Proposition 3.2.9 in \onlinecite{Schaumann2013} where a proof is given. Proposition~\ref{eno_remark_3.9} also appears as Remark~3.9 in \onlinecite{MR2677836} where it is used to give a formula for the Brauer-Picard ring of $\vvec{A}$ when $A$ is an abelian group. The other ingredient in the proof is the classical Morita's theorem (Proposition~\ref{prop:cmt}). A proof is given in \onlinecite{1610.05204}.
\begin{proposition}\label{prop:cmt}
Let $\mathcal{A}$ be an abelian category and $\mathcal{P} \subseteq \mathcal{A}$ a full subcategory consisting of projective objects that generate $\mathcal{A}$. Then $\mathcal{A} \cong [\mathcal{P}^{\rm op}, \vvec{}]$.
\end{proposition}
{\em Proof of Proposition~\ref{ladder_model_relative_tensor_product}: }
Let $\mathcal{B}$ be a fusion category and $\mathcal{M} \curvearrowleft \mathcal{B} \curvearrowright \mathcal{N}$ be bimodules. Proposition~\ref{eno_remark_3.9} implies that the $I$-modules $I \otimes m \boxtimes n$ are a generating set of projectives in the relative tensor product $\mathcal{M} \boxtimes_{\mathcal{B}} \mathcal{N}$. Define $\lad{\mathcal{M},\mathcal{N}}{\mathcal{B}}$ to be the full subcategory of $\mathcal{M} \boxtimes_{\mathcal{B}} \mathcal{N}$ with objects $I \otimes m \boxtimes n$ where $m$ varies over a complete set of simples in $\mathcal{M}$ and $n$ varies over a complete set of simples in $\mathcal{N}$. By Morita's theorem, $\mathcal{M} \boxtimes_{\mathcal{B}} \mathcal{N}$ is equivalent to the category of representations of $\lad{\mathcal{M},\mathcal{N}}{\mathcal{B}}^{\rm op}$. To produce a nice model $\lad{\mathcal{M},\mathcal{N}}{\mathcal{B}}$, we need to understand the multiplication map on $I$. We have
\begin{align}I^{\otimes 2} = \bigoplus_{b,c \in \Lambda} c^* b^* \boxtimes b c.\end{align}
Post composing $c^* b^* \boxtimes b c \to I^{\otimes 2}$ with the multiplication map gives the vector corresponding to decomposing the identity on $bc$ in
\begin{align} \mathcal{B}^{\otimes {\rm op}} \boxtimes \mathcal{B} \left( c^*b^* \boxtimes bc, \bigoplus_{a \in \Lambda} a^* \boxtimes a \right) = \bigoplus_{a \in \Lambda} \mathcal{B}(a,bc) \otimes \mathcal{B}(bc,a).\end{align}
By definition, we have
\begin{align} {\bf Mod}(I)(I \otimes m \boxtimes n, I \otimes x \boxtimes y) &= \mathcal{M} \boxtimes \mathcal{N}(m \boxtimes n, I \otimes x \boxtimes y)
  \\
  &= \bigoplus_{b \in \Lambda} \mathcal{M} \boxtimes \mathcal{N}(m \boxtimes n,x b^* \boxtimes by) \\
  &= \bigoplus_{b \in \Lambda} \mathcal{M}(m,xb^*) \otimes \mathcal{N}(n,by) \\
  &= \bigoplus_{b \in \Lambda} \mathcal{M}(mb,x) \otimes \mathcal{N}(n,by)
\end{align}
We can depicted rank 1 tensors in $\mathcal{M}(mb,x) \otimes \mathcal{N}(n,by)$ using {\bf ladder} string diagrams of the following form
\begin{align}
\begin{array}{c}
\includeTikz{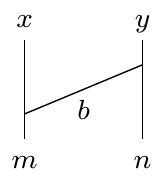}
{
	\begin{tikzpicture}
	\coordinate (A) at (0,0);
	\ladder[A][$m$][$n$][$b$][$x$][$y$];
	\end{tikzpicture}
}
\end{array}.
\end{align}
Composition in $\lad{\mathcal{M},\mathcal{N}}{\mathcal{B}}$ corresponds to stacking the ladders and then decomposing the rungs into H-diagrams. \qed

The computations in Section~\ref{S:Zp} take the following form: First, we compute the ladder category for a given pair of bimodules $\mathcal{M}$ and $\mathcal{N}$. Next, the isomorphism classes of simple objects in the Karoubi envelope of $\lad{\mathcal{M},\mathcal{N}}{\vvec{\ZZ{p}}}$ are identified. Finally, we compute the left and right actions of $\mathcal{A}$ and $\mathcal{C}$ respectively, and the associator. This allows the corresponding indecomposable bimodules to be identified.

The left action of $\mathcal{A}$ on the ladder category is defined by
\begin{align}
a\otimes
\begin{array}{c}
\includeTikz{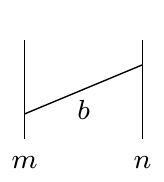}
{
	\begin{tikzpicture}
		\coordinate (A) at (0,0);
		\ladder[A][$m$][$n$][$b$][][];
	\end{tikzpicture}
}
\end{array}
=
\begin{array}{c}
\includeTikz{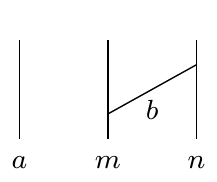}
{
	\begin{tikzpicture}[xscale=.75]
	\coordinate (A) at (0,0);
	\Lladder[A][$m$][$n$][$b$][][][$a$];
	\end{tikzpicture}
}
\end{array}.
\end{align}
The right action is defined using a similar diagram. There is a slight subtlety related to the left and right actions of $\mathcal{A}$ and $\mathcal{C}$ respectively. Indeed, the actions of $\mathcal{A}$ and $\mathcal{C}$ do not actually extend from $\mathcal{M}$ and $\mathcal{N}$ to the ladder category $\lad{\mathcal{M},\mathcal{N}}{\vvec{\ZZ{p}}}$ as we have defined it. This problem is easy to fix. Instead of defining the objects in the ladder category to be pairs of simples, we should really define the objects to be pairs of arbitrary tensor products of simples. Making this change adds additional complexity to the definitions and it doesn't change any of the computations as we will see below.

Once the ladder category has been suitably described, we need to compute its Karoubi envelope. The first step is to find a representative for each simple object in the envelope. If $x$ is such a representative, and $a \in \mathcal{A}$ is simple, then we need to find a basis for all the morphisms between $a \otimes x$ and the other representatives in the envelope. These morphisms will be our trivalent vertices in the Karoubi envelope and using them we can identify the fusion rules and the associators. 

%

\section{Computing the Brauer-Picard Ring of $\vvec{\ZZ{p}}$}\label{S:Zp}

\begin{table}[t]
	\resizebox{\linewidth}{!}{
		\begin{tabular}{!{\vrule width 1pt}>{\columncolor[gray]{.9}[\tabcolsep]}c!{\vrule width 1pt}c !{\color[gray]{.8}\vrule} c !{\color[gray]{.8}\vrule} c !{\color[gray]{.8}\vrule} c !{\color[gray]{.8}\vrule} c !{\color[gray]{.8}\vrule} c!{\vrule width 1pt}}
			\toprule[1pt]
			\rowcolor[gray]{.9}[\tabcolsep]Bimodule label & Subgroup & \# objects & Cosets $\to$ object label & Left action & Right action & Associator \\
			\toprule[1pt]
			$T$ & $\{(0,0)\}$ & $p^2$ & $\{(g,h)\}\to(g,h)$&
			$
			\begin{array}{c}
			\includeTikz{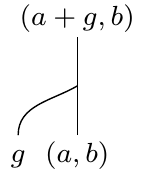}
			{
				\begin{tikzpicture}
				\Laction[$g$][$(a,b)$][$(a+g,b)$];
				\end{tikzpicture}
			}
			\end{array}
			$
			&
			$
			\begin{array}{c}
			\includeTikz{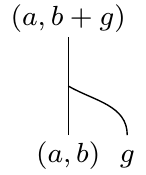}
			{
				\begin{tikzpicture}
				\Raction[$g$][$(a,b)$][$(a,b+g)$];
				\end{tikzpicture}
			}
			\end{array}
			$&
			$
			\begin{array}{c}
			\includeTikz{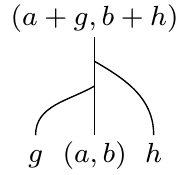}
			{
				\begin{tikzpicture}
				\Lassociator[$g$][$h$][$(a,b)$][$(a+g,b+h)$];
				\end{tikzpicture}
			}
			\end{array}
			=
			\begin{array}{c}
			\includeTikz{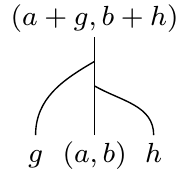}
			{
				\begin{tikzpicture}
				\Rassociator[$g$][$h$][$(a,b)$][$(a+g,b+h)$];
				\end{tikzpicture}
			}
			\end{array}
			$
			\\
			\greycline{2-7}
			$L$ & $\langle (1,0) \rangle\cong\ZZ{p}$ & p & $\{(h,g)\,|\,h\in \ZZ{p}\}\to g$&
			$
			\begin{array}{c}
			\includeTikz{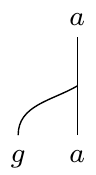}
			{
				\begin{tikzpicture}
				\Laction[$g$][$a$][$a$];
				\end{tikzpicture}
			}
			\end{array}
			$
			&
			$
			\begin{array}{c}
			\includeTikz{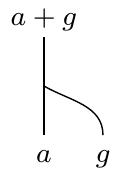}
			{
				\begin{tikzpicture}
				\Raction[$g$][$a$][$a+g$];
				\end{tikzpicture}
			}
			\end{array}
			$&
			$
			\begin{array}{c}
			\includeTikz{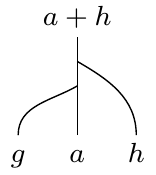}
			{
				\begin{tikzpicture}
				\Lassociator[$g$][$h$][$a$][$a+h$];
				\end{tikzpicture}
			}
			\end{array}
			=
			\begin{array}{c}
			\includeTikz{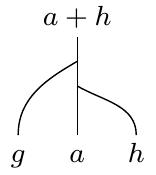}
			{
				\begin{tikzpicture}
				\Rassociator[$g$][$h$][$a$][$a+h$];
				\end{tikzpicture}
			}
			\end{array}
			$
			\\
			\greycline{2-7} 
			$R$ & $\langle (0,1) \rangle\cong\ZZ{p}$ & p & $\{(g,h)\,|\,h\in \ZZ{p}\}\to g$&
			$
			\begin{array}{c}
			\includeTikz{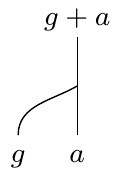}
			{
				\begin{tikzpicture}
				\Laction[$g$][$a$][$g+a$];
				\end{tikzpicture}
			}
			\end{array}
			$
			&
			$
			\begin{array}{c}
			\includeTikz{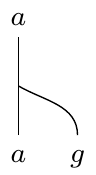}
			{
				\begin{tikzpicture}
				\Raction[$g$][$a$][$a$];
				\end{tikzpicture}
			}
			\end{array}
			$&
			$
			\begin{array}{c}
			\includeTikz{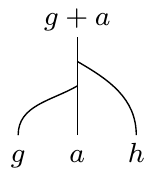}
			{
				\begin{tikzpicture}
				\Lassociator[$g$][$h$][$a$][$g+a$];
				\end{tikzpicture}
			}
			\end{array}
			=
			\begin{array}{c}
			\includeTikz{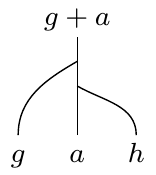}
			{
				\begin{tikzpicture}
				\Rassociator[$g$][$h$][$a$][$g+a$];
				\end{tikzpicture}
			}
			\end{array}
			$
			\\
			\greycline{2-7} 
			$F_0$ & $\langle (0,1), (1,0) \rangle\cong\ZZ{p}\times\ZZ{p}$ & 1 & $\ZZ{p}\times \ZZ{p}\to *$&
			$
			\begin{array}{c}
			\includeTikz{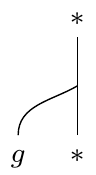}
			{
				\begin{tikzpicture}
				\Laction[$g$][$*$][$*$];
				\end{tikzpicture}
			}
			\end{array}
			$
			&
			$
			\begin{array}{c}
			\includeTikz{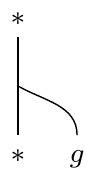}
			{
				\begin{tikzpicture}
				\Raction[$g$][$*$][$*$];
				\end{tikzpicture}
			}
			\end{array}
			$&
			$
			\begin{array}{c}
			\includeTikz{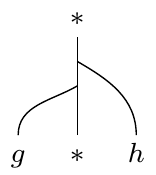}
			{
				\begin{tikzpicture}
				\Lassociator[$g$][$h$][$*$][$*$];
				\end{tikzpicture}
			}
			\end{array}
			=
			\begin{array}{c}
			\includeTikz{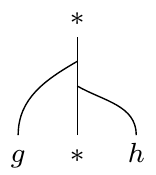}
			{
				\begin{tikzpicture}
				\Rassociator[$g$][$h$][$*$][$*$];
				\end{tikzpicture}
			}
			\end{array}
			$
			\\		
			\toprule[1pt]
			$X_k$ & $\langle (-k,1) \rangle\cong\ZZ{p}$ & p & $\{n(-k,1)+(h,0)\,|\,n\in \ZZ{p}\}\to h$&
			$
			\begin{array}{c}
			\includeTikz{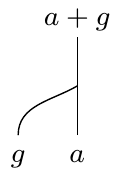}
			{
				\begin{tikzpicture}
				\Laction[$g$][$a$][$a+g$];
				\end{tikzpicture}
			}
			\end{array}
			$
			&
			$
			\begin{array}{c}
			\includeTikz{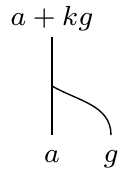}
			{
				\begin{tikzpicture}
				\Raction[$g$][$a$][$a+kg$];
				\end{tikzpicture}
			}
			\end{array}
			$&
			$
			\begin{array}{c}
			\includeTikz{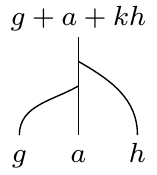}
			{
				\begin{tikzpicture}
				\Lassociator[$g$][$h$][$a$][$g+a+kh$];
				\end{tikzpicture}
			}
			\end{array}
			=
			\begin{array}{c}
			\includeTikz{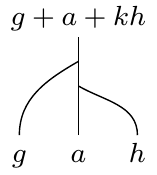}
			{
				\begin{tikzpicture}
				\Rassociator[$g$][$h$][$a$][$g+a+kh$];
				\end{tikzpicture}
			}
			\end{array}
			$
			\\
			\greycline{2-7}
			$F_q$& $\langle (0,1), (1,0) \rangle\cong\ZZ{p}\times\ZZ{p}$ & 1 & $\ZZ{p}\times \ZZ{p}\to *$&
			$
			\begin{array}{c}
			\includeTikz{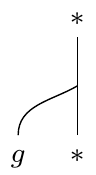}
			{
				\begin{tikzpicture}
				\Laction[$g$][$*$][$*$];
				\end{tikzpicture}
			}
			\end{array}
			$
			&
			$
			\begin{array}{c}
			\includeTikz{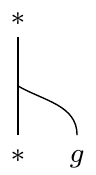}
			{
				\begin{tikzpicture}
				\Raction[$g$][$*$][$*$];
				\end{tikzpicture}
			}
			\end{array}
			$&
			$
			\begin{array}{c}
			\includeTikz{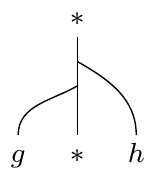}
			{
				\begin{tikzpicture}
				\Lassociator[$g$][$h$][$*$][$*$];
				\end{tikzpicture}
			}
			\end{array}
			=
			e^{\frac{2\pi i}{p} q g h}
			\begin{array}{c}
			\includeTikz{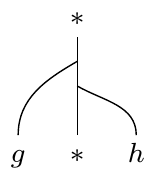}
			{
				\begin{tikzpicture}
				\Rassociator[$g$][$h$][$*$][$*$];
				\end{tikzpicture}
			}
			\end{array}
			$
			\\
			\hline
		\end{tabular}
	}
	\caption{Data for all $\vvec{\ZZ{p}}-\vvec{\ZZ{p}}$ bimodules. $q \in H^2(\ZZ{p},U(1)) \cong \ZZ{p}$. Bimodules below the thick line are invertible.}\label{tab:zpdata}
\end{table}

In this section, we compute the Brauer-Picard ring for $\vvec{\ZZ{p}}$ with $p$ prime. We remark that the ring for the twisted case $\vvectwist{\ZZ{p}}{\phi}$ can be similarly computed. In that case, there are additional consistency conditions between the associators and cocycle twist that restrict the bimodules compared to $\vvec{\ZZ{p}}$.

Recall from Chapter~7 of \onlinecite{MR3242743} that each $\vvec{G}\text{ -}\vvec{G}$ bimodule $\mathcal{M}$ is labelled by a conjugacy class of subgroups $H\subseteq G\times G$ and a 2-cocycle $\omega\in H^2(H,U(1))$. The simple objects in $\mathcal{M}$ are labelled by cosets of $H$. In Table~\ref{tab:zpdata}, we list the bimodules, and the associated functors ($\vartriangleright,\, \vartriangleleft$ in Definition~\ref{bimodule_definition}). The associators are obtained from the cocycle $\omega$ using the formula explained in \onlinecite{288303}. We now provide some example calculations. The full multiplication table is presented in Table~\ref{tab:zptable}.

\begin{exmp}[$T \otimes_{\vvec{\ZZ{p}}} T$]
 We begin by computing $T\otimes_{\vvec{\ZZ{p}}} T$. We proceed in two steps: First we construct the \emph{ladder category} $\lad{T,T}{}$, then we construct its \emph{Karoubi envelope} $\kar{\protect\lad{T,T}{}}$ and identify the left and right actions and associators. With this data, we can identify the resulting bimodule. This example is unnaturally simple because the ladder category is already idempotent complete. The (simple) objects in $\lad{T, T}{}$ are given by pairs of objects from $T$, written $(g_0,g_1)(g_2,g_2)$ with $g_i\in \ZZ{p}$. Therefore, $\lad{T,T}{}$ has $p^4$ simple objects. Given two objects $(a,b)(c,d)$ and $(a,b-g)(c+g,d)$, the full set of morphisms between them are
\begin{align}
\begin{array}{c}
\includeTikz{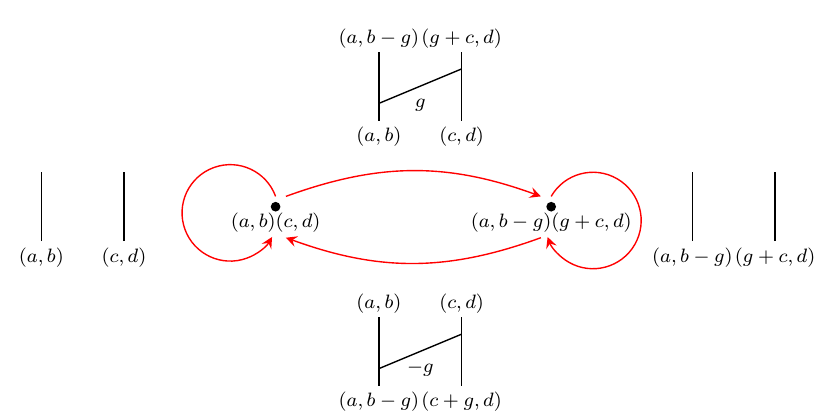}{
	\begin{tikzpicture}[scale=0.7, every node/.style={scale=0.7}]
	\clip (-4,-3) rectangle (8,3);
	\coordinate (A) at (0,0);
	\coordinate (B) at (4,0);
	\coordinate (C) at (-4,0);
	\coordinate (D) at (3,2.5);
	\coordinate (E) at (9.5,0);
	\coordinate (F) at (3,-3);
	\fill (A) circle(2pt) node[below,inline text] {\footnotesize$(a,b)(c,d)$};
	\fill (B) circle(2pt) node[below,inline text] {\footnotesize$(a,b-g)(g+c,d)$};
	\draw[-stealth,red] ($(A)+(0,.15)$) arc (20:330:.7);
	\draw[-stealth,red] ($(B)+(0,.15)$) arc (330:20:-.7);
	\draw[-stealth,red] ($(A)+(.15,.15)$) to[out=20,in=180-20] ($(B)+(-.15,.15)$);
	\draw[-stealth,red] ($(B)+(-.15,-.45)$) to[out=180+20,in=-20] ($(A)+(.15,-.45)$);
	\ladder[C][$(a,b)$][$(c,d)$][][][];
	\ladder[E][$(a,b-g)$][$(g+c,d)$][][][];
	\ladder[D][$(a,b)$][$(c,d)$][$g$][$(a,b-g)$][$(g+c,d)$];
	\ladder[F][$(a,b-g)$][$(c+g,d)$][$-g$][$(a,b)$][$(c,d)$];
	\end{tikzpicture}
}
\end{array}.\label{eq:TTmorph}
\end{align}
For all other pairs of objects, only the self-morphisms exist. Since all morphisms in Equation~\ref{eq:TTmorph} are invertible, the objects $(a,b)(c,d)$ and $(a,b-g)(c+g,d)$ are isomorphic in $\lad{T,T}{}$. The objects in the Karoubi envelope of the ladder category are labelled by pairs $(A,e)$, where $A$ is an object from $\lad{T,T}{}$ and $e:A\to A$ and idempotent. The $p^3$ non-isomorphic objects of the Karoubi envelope can be represented by
\begin{align}
\begin{array}{c}
\includeTikz{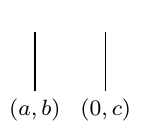}
{
	\begin{tikzpicture}[scale=0.6, every node/.style={scale=0.8}]
	\coordinate (A) at (0,0);
	\ladder [A][$(a,b)$][$(0,c)$][][][];
	\end{tikzpicture}
}
\end{array}.
\end{align}

\begin{table}[t]
	\begin{tabular}{!{\vrule width 1pt}>{\columncolor[gray]{.9}[\tabcolsep]}c!{\vrule width 1pt}c c c c!{\vrule width 1pt}c c!{\vrule width 1pt}}
		\toprule[1pt]
		\rowcolor[gray]{.9}[\tabcolsep]$\otimes_{\vvec{\ZZ{p}}}$&$T$&$L$&$R$&$F_0$&$X_l$&$F_r$\\
		\toprule[1pt]
		$T$&$p\cdot T$&$T$&$p\cdot R$&$R$&$T$&$R$\\
		$L$&$p\cdot L$&$L$&$p\cdot F_0$&$F_0$&$L$&$F_0$\\
		$R$&$T$&$p\cdot T$&$R$&$p\cdot R$&$R$&$T$\\
		$F_0$&$L$&$p\cdot L$&$F_0$&$p\cdot F_0$&$F_0$&$L$\\
		\toprule[1pt]
		$X_k$&$T$&$L$&$R$&$F_0$&$X_{kl}$&$F_{k^{-1}r}$\\
		$F_q$&$L$&$T$&$F_0$&$R$&$F_{ql}$&$X_{q^{-1}r}$\\
		\toprule[1pt]
	\end{tabular}
	\caption{Multiplication table for $\protect\bpr{\protect\vvec{\ZZ{p}}}$.}\label{tab:zptable}
\end{table}

We now need to identify the left action of $\vvec{\ZZ{p}}$. A morphism $f:(A,e)\to(A^\prime,e^\prime)$ is a morphism $f\in \lad{T,T}{}$ which obeys $ef=f=fe^\prime$. Such a morphism allows us to identify the left action of $g$
\begin{align}
\begin{array}{c}
\includeTikz{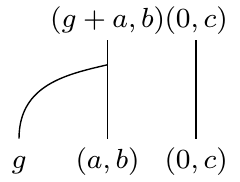}
{
	\begin{tikzpicture}[xscale=.75]
	\coordinate (A) at (0,0);
	\Lmorphism[A][$(a,b)$][$(0,c)$][][$(g+a,b)$][$(0,c)$][$g$];
	\end{tikzpicture}
}
\end{array}
:
\begin{array}{c}
\includeTikz{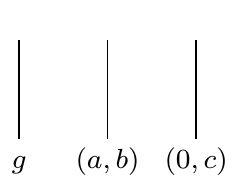}
{
	\begin{tikzpicture}[xscale=.75]
	\coordinate (A) at (0,0);
	\Lladder[A][$(a,b)$][$(0,c)$][][][][$g$];
	\end{tikzpicture}
}
\end{array}
\tilde{\to}
\begin{array}{c}
\includeTikz{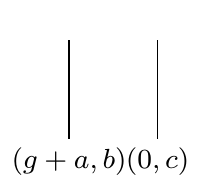}
{
	\begin{tikzpicture}[xscale=.75]
	\coordinate (A) at (0,0);
	\ladder[A][$(g+a,b)$][$(0,c)$][][][];
	\end{tikzpicture}
}
\end{array}.\label{eqn:TTleftiso}
\end{align}
One can readily check that this morphism ``absorbs'' the appropriate idempotents and therefore establishes an isomorphism. The left action is therefore $g \otimes (a,b)(0,c)=(g+a,b)(0,c)$. A similar isomorphism 
\begin{align}
\begin{array}{c}
\includeTikz{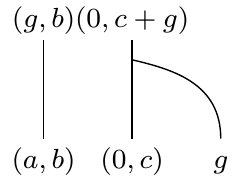}
{
	\begin{tikzpicture}[xscale=.75]
	\coordinate (A) at (0,0);
	\Rmorphism[A][$(a,b)$][$(0,c)$][][$(g,b)$][$(0,c+g)$][$g$];
	\end{tikzpicture}
}
\end{array}
:
\begin{array}{c}
\includeTikz{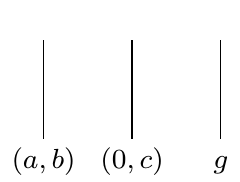}
{
	\begin{tikzpicture}[xscale=.75]
	\coordinate (A) at (0,0);
	\Rladder[A][$(a,b)$][$(0,c)$][][][][$g$];
	\end{tikzpicture}
}
\end{array}
\tilde{\to}
\begin{array}{c}
\includeTikz{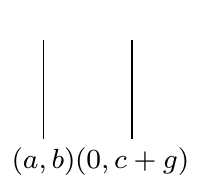}
{
	\begin{tikzpicture}[xscale=.75]
	\coordinate (A) at (0,0);
	\ladder[A][$(a,b)$][$(0,c+g)$][][][];
	\end{tikzpicture}
}
\end{array},\label{eqn:TTrightiso}
\end{align}
which identifies the right action as $(a,b)(0,c)\otimes g=(a,b)(0,c+g)$. These actions are consistent with the bimodule $T$. Since there is no way to change the label $b$, it labels $p$ different blocks with no morphisms between them. It is straightforward to check that the morphisms in Equation~\ref{eqn:TTleftiso} and Equation~\ref{eqn:TTrightiso} commute, meaning that the associator is trivial. It follows that
\begin{align}
	T\otimes_{\vvec{\ZZ{p}}}T =p T.
\end{align}
\end{exmp}
\begin{center}\rule{.8\textwidth}{1pt}\end{center}

\begin{exmp}[$R \otimes_{\vvec{\ZZ{p}}} F_0$]
The simple objects in $\lad{R,F_0}{}$ are labelled by $a$, with $a\in \ZZ{p}$, so there are $p$ simple objects. Given an object $a\in \lad{R,F_0}{}$, there is a $p$ dimensional algebra of morphisms
\begin{align}
\begin{array}{c}
\includeTikz{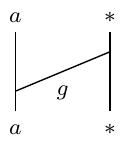}
{
	\begin{tikzpicture}[scale=0.8, every node/.style={scale=0.8}]
	\coordinate (A) at (0,0);
	\ladder[A][$a$][$*$][$g$][$a$][$*$];
	\end{tikzpicture}
}
\end{array}&,\\
\text{with }
\begin{array}{c}
\includeTikz{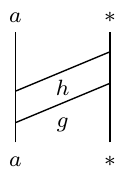}
{
	\begin{tikzpicture}[scale=0.8, every node/.style={scale=0.8}]
	\coordinate (A) at (0,0);\coordinate (B) at (0,.5);
	\ladder[A][$a$][$*$][$g$][][];
	\ladder[B][][][$h$][$a$][$*$];
	\end{tikzpicture}
}
\end{array}
&=\begin{array}{c}
\includeTikz{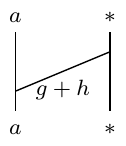}
{
	\begin{tikzpicture}[scale=0.8, every node/.style={scale=0.8}]
	\coordinate (A) at (0,0);
	\ladder[A][$a$][$*$][$g+h$][$a$][$*$];
	\end{tikzpicture}
}
\end{array}
\end{align}
which form the group algebra $\mathbb{C}\left[\ZZ{p}\right]$. The $p$ indecomposable idempotents $a\to a$ are given by
\begin{align}
I_{k}&=\frac{1}{p}\sum_{g\in\ZZ{p}}
e^{\frac{2\pi i}{p}kg}
\begin{array}{c}
\includeTikz{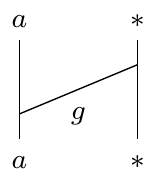}
{
	\begin{tikzpicture}
	\coordinate (A) at (0,0);
	\ladder[A][$a$][$*$][$g$][$a$][$*$];
	\end{tikzpicture}
}
\end{array}.
\end{align}
The $p^2$ simple objects of $\kar{\lad{R,F_0}{}}$ can be labelled by $(a,k)$. To compute the action of $g$ on these objects, we again need to find a morphism $f:g\otimes(a,k)\to(a^\prime,k^\prime)$ in $\kar{\lad{R,F_0}{}}$. In this case, we can take $k' = k$ and $a' = g + a$. Recall that this needs to absorb the idempotents $I_k$ and $I_k^\prime$ on the bottom and top respectively. Such a morphism is given by 
\begin{align}
\frac{1}{p}\sum_{x\in\ZZ{p}}e^{\frac{2\pi i}{p}kx}
\begin{array}{c}
\includeTikz{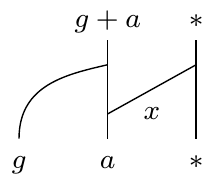}
{
	\begin{tikzpicture}[xscale=.75]
	\coordinate (A) at (0,0);
	\Lmorphism[A][$a$][$*$][$x$][$g+a$][$*$][$g$];
	\end{tikzpicture}
}
\end{array}
:
\frac{1}{p}\sum_{y\in\ZZ{p}}e^{\frac{2\pi i}{p}ky}
\begin{array}{c}
\includeTikz{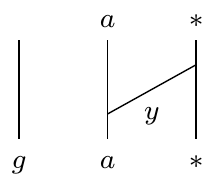}
{
	\begin{tikzpicture}[xscale=.75]
	\coordinate (A) at (0,0);
	\Lladder[A][$a$][$*$][$y$][$a$][$*$][$g$];
	\end{tikzpicture}
}
\end{array}
\tilde{\to}
\frac{1}{p}\sum_{y\in\ZZ{p}}e^{\frac{2\pi i}{p}ky}
\begin{array}{c}
\includeTikz{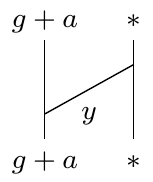}
{
	\begin{tikzpicture}[xscale=.75]
	\coordinate (A) at (0,0);
	\ladder[A][$g+a$][$*$][$y$][$g+a$][$*$];
	\end{tikzpicture}
}
\end{array}.\label{eqn:RF0leftiso}
\end{align}
One can check that the idempotents are absorbed:
\begin{align}
	I_kf&=
	\frac{1}{p^2}
	\sum_{x,y\in \ZZ{p}}
	e^{\frac{2\pi i}{p}k(x+y)}
	\begin{array}{c}
		\includeTikz{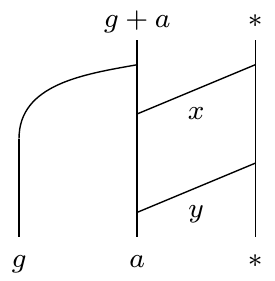}
		{
			\begin{tikzpicture}
				\coordinate (A) at (0,0);
				\coordinate (B) at (0,-1);
				\Lmorphism[A][][][$x$][$g+a$][$*$][];
				\Lladder[B][$a$][$*$][$y$][][][$g$];
			\end{tikzpicture}
		}
	\end{array}
	=
	\frac{1}{p^2}
	\sum_{x,y\in \ZZ{p}}
	e^{\frac{2\pi i}{p}k(x+y)}
	\begin{array}{c}
	\includeTikz{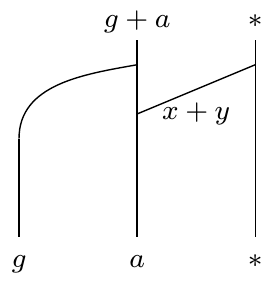}
	{
		\begin{tikzpicture}
		\coordinate (A) at (0,0);
		\coordinate (B) at (0,-1);
		\Lmorphism[A][][][$x+y$][$g+a$][$*$][];
		\Lladder[B][$a$][$*$][][][][$g$];
		\end{tikzpicture}
	}
	\end{array}\\
	&=
	\frac{1}{p^2}
	\sum_{x,z\in \ZZ{p}}
	e^{\frac{2\pi i}{p}kz}
	\begin{array}{c}
	\includeTikz{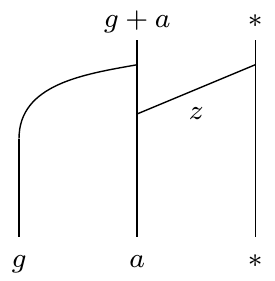}
	{
		\begin{tikzpicture}
		\coordinate (A) at (0,0);
		\coordinate (B) at (0,-1);
		\Lmorphism[A][][][$z$][$g+a$][$*$][];
		\Lladder[B][$a$][$*$][][][][$g$];
		\end{tikzpicture}
	}
	\end{array}
	=
	\frac{1}{p}
	\sum_{z\in \ZZ{p}}
	e^{\frac{2\pi i}{p}kz}
	\begin{array}{c}
	\includeTikz{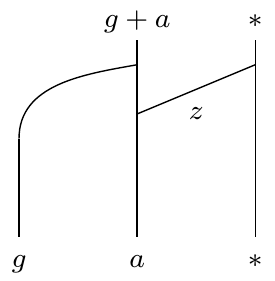}
	{
		\begin{tikzpicture}
		\coordinate (A) at (0,0);
		\coordinate (B) at (0,-1);
		\Lmorphism[A][][][$z$][$g+a$][$*$][];
		\Lladder[B][$a$][$*$][][][][$g$];
		\end{tikzpicture}
	}
	\end{array}
	\\
	&=f.
\end{align}
Absorption from above can be checked similarly, recalling that the associator on the left strand is trivial
\begin{align}
	fI_{k+g}&=
	\frac{1}{p^2}
	\sum_{x,y\in \ZZ{p}}
	e^{\frac{2\pi i}{p}k(x+y)}
	\begin{array}{c}
		\includeTikz{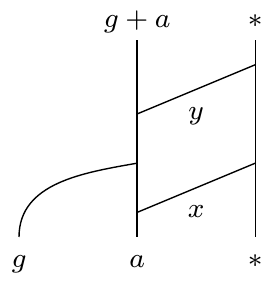}
		{
			\begin{tikzpicture}
				\coordinate (A) at (0,0);
				\coordinate (B) at (0,-1);
				\Lmorphism[B][$a$][$*$][$x$][][][$g$];
				\ladder[A][][][$y$][$g+a$][$*$];
			\end{tikzpicture}
		}
	\end{array}
	=
	\frac{1}{p^2}
	\sum_{x,y\in \ZZ{p}}
	e^{\frac{2\pi i}{p}k(x+y)}
	\begin{array}{c}
		\includeTikz{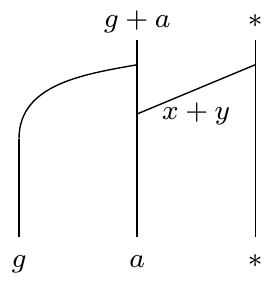}
		{
			\begin{tikzpicture}
			\coordinate (A) at (0,0);
			\coordinate (B) at (0,-1);
			\Lmorphism[A][][][$x+y$][$g+a$][$*$][];
			\Lladder[B][$a$][$*$][][][][$g$];
			\end{tikzpicture}
		}
	\end{array}
	=f.
\end{align}

For the right action, we find
\begin{align}
\frac{1}{p}\sum_{x\in\ZZ{p}}e^{\frac{2\pi i}{p}kx}
\begin{array}{c}
\includeTikz{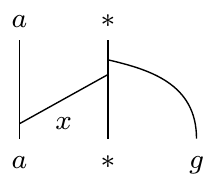}
{
	\begin{tikzpicture}[xscale=.75]
	\coordinate (A) at (0,0);
	\Rmorphism[A][$a$][$*$][$x$][$a$][$*$][$g$];
	\end{tikzpicture}
}
\end{array}
:
\frac{1}{p}\sum_{y\in\ZZ{p}}e^{\frac{2\pi i}{p}ky}
\begin{array}{c}
\includeTikz{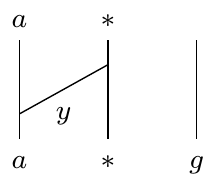}
{
	\begin{tikzpicture}[xscale=.75]
	\coordinate (A) at (0,0);
	\Rladder[A][$a$][$*$][$y$][$a$][$*$][$g$];
	\end{tikzpicture}
}
\end{array}
\tilde{\to}
\frac{1}{p}\sum_{y\in\ZZ{p}}e^{\frac{2\pi i}{p}ky}
\begin{array}{c}
\includeTikz{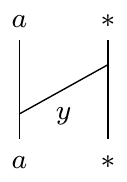}
{
	\begin{tikzpicture}[xscale=.75]
	\coordinate (A) at (0,0);
	\ladder[A][$a$][$*$][$y$][$a$][$*$];
	\end{tikzpicture}
}
\end{array}.\label{eqn:RF0rightiso}
\end{align}
The left action is therefore $g\otimes (a,k)=(g+a,k)$, with trivial right action. As before, the morphisms Equation~\ref{eqn:RF0leftiso} and Equation.~\ref{eqn:RF0rightiso} commute so the associator is trivial. The label $k$ is unaffected by the group action, so labels $p$ distinct blocks with no morphisms between them. This gives
\begin{align}
R\otimes_{\vvec{\ZZ{p}}} F_0=pR.
\end{align}

\end{exmp}
\begin{center}\rule{.8\textwidth}{1pt}\end{center}

\begin{exmp}[$X_k\otimes_{\vvec{\ZZ{p}}} X_l$]
The simple objects in $\lad{X_k,X_l}{}$ are labeled by $(a,b)$, with $a,b\in \ZZ{p}$. Given a pair of objects $(a,b)$ and $(a-kg,b+g)$ the full set of morphisms between them are
\begin{align}
	\begin{array}{c}
		\includeTikz{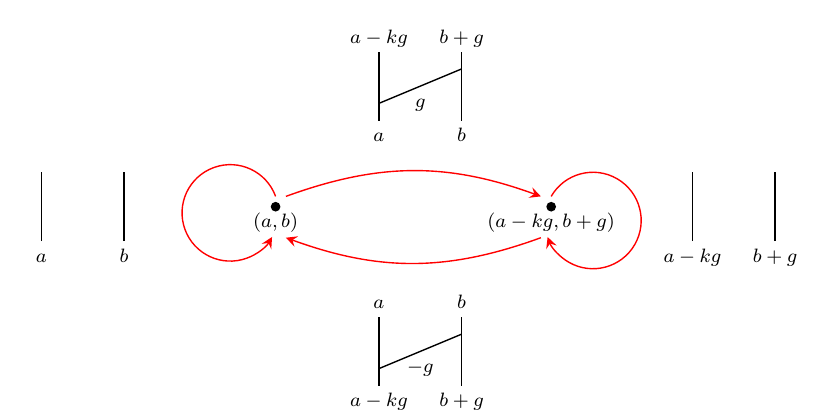}{
			\begin{tikzpicture}[scale=0.7, every node/.style={scale=0.7}]
			\clip (-4,-3) rectangle (8,3);
			\coordinate (A) at (0,0);
			\coordinate (B) at (4,0);
			\coordinate (C) at (-4,0);
			\coordinate (D) at (3,2.5);
			\coordinate (E) at (9.5,0);
			\coordinate (F) at (3,-3);
			\fill (A) circle(2pt) node[below,inline text] {\footnotesize$(a,b)$};
			\fill (B) circle(2pt) node[below,inline text] {\footnotesize$(a-kg,b+g)$};
			\draw[-stealth,red] ($(A)+(0,.15)$) arc (20:330:.7);
			\draw[-stealth,red] ($(B)+(0,.15)$) arc (330:20:-.7);
			\draw[-stealth,red] ($(A)+(.15,.15)$) to[out=20,in=180-20] ($(B)+(-.15,.15)$);
			\draw[-stealth,red] ($(B)+(-.15,-.45)$) to[out=180+20,in=-20] ($(A)+(.15,-.45)$);
			\ladder[C][$a$][$b$][][][];
			\ladder[E][$a-kg$][$b+g$][][][];
			\ladder[D][$a$][$b$][$g$][$a-kg$][$b+g$];
			\ladder[F][$a-kg$][$b+g$][$-g$][$a$][$b$];
			\end{tikzpicture}
		}
	\end{array}.\label{eq:XXmorph}
\end{align}
For all other pairs of objects, only the self-morphisms exist. Since all morphisms in Equation~\ref{eq:XXmorph} are invertible, the objects $(a,b)$ and $(a-kg,b+g)$ are isomorphic. The $p$ non-isomorphic objects of $\kar{\lad{X_k,X_l}{}}$ can be represented by
\begin{align}
	\begin{array}{c}
		\includeTikz{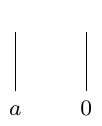}
		{
			\begin{tikzpicture}[scale=0.6, every node/.style={scale=0.8}]
			\coordinate (A) at (0,0);
			\ladder [A][$a$][$0$][][][];
			\end{tikzpicture}
		}
	\end{array}.
\end{align}
The morphism associated to the left $g$ action is
\begin{align}
	\begin{array}{c}
		\includeTikz{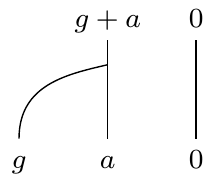}
		{
			\begin{tikzpicture}[xscale=.75]
			\coordinate (A) at (0,0);
			\Lmorphism[A][$a$][$0$][][$g+a$][$0$][$g$];
			\end{tikzpicture}
		}
	\end{array}
	:
	\begin{array}{c}
		\includeTikz{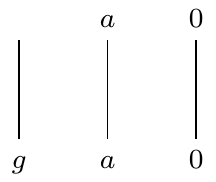}
		{
			\begin{tikzpicture}[xscale=.75]
			\coordinate (A) at (0,0);
			\Lladder[A][$a$][$0$][][$a$][$0$][$g$];
			\end{tikzpicture}
		}
	\end{array}
	\tilde{\to}
	\begin{array}{c}
		\includeTikz{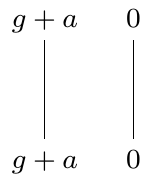}
		{
			\begin{tikzpicture}[xscale=.75]
			\coordinate (A) at (0,0);
			\ladder[A][$g+a$][$0$][][$g+a$][$0$];
			\end{tikzpicture}
		}
	\end{array},\label{eqn:XXleftiso}
\end{align}
so $g\otimes(a,0)=(a+g,0)$. The right action is more complicated. Consider the morphism
\begin{align}
	\begin{array}{c}
		\includeTikz{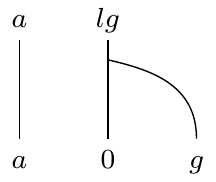}
		{
			\begin{tikzpicture}[xscale=.75]
			\coordinate (A) at (0,0);
			\Rmorphism[A][$a$][$0$][][$a$][$lg$][$g$];
			\end{tikzpicture}
		}
	\end{array}
	:
	\begin{array}{c}
		\includeTikz{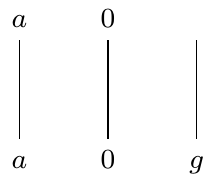}
		{
			\begin{tikzpicture}[xscale=.75]
			\coordinate (A) at (0,0);
			\Rladder[A][$a$][$0$][][$a$][$0$][$g$];
			\end{tikzpicture}
		}
	\end{array}
	\tilde{\to}
	\begin{array}{c}
		\includeTikz{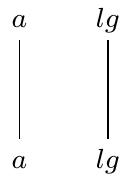}
		{
			\begin{tikzpicture}[xscale=.75]
			\coordinate (A) at (0,0);
			\ladder[A][$a$][$lg$][][$a$][$lg$];
			\end{tikzpicture}
		}
	\end{array}.\label{eqn:RF0rightiso1}
\end{align}
The codomain is not the chosen representative of the isomorphism class. Instead, we need to append the morphism
\begin{align}
	\begin{array}{c}
		\includeTikz{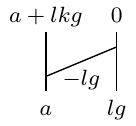}
		{
			\begin{tikzpicture}[scale=0.6, every node/.style={scale=0.8}]
			\coordinate (A) at (0,0);
			\ladder [A][$a$][$lg$][$-lg$][$a+lkg$][$0$];
			\end{tikzpicture}
		}
	\end{array}
\end{align}
to identify the correct block. The right action is then identified from
\begin{align}
	\begin{array}{c}
		\includeTikz{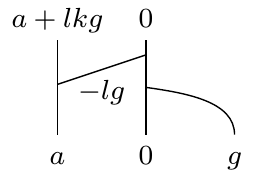}
		{
			\begin{tikzpicture}[xscale=.75,yscale=.6]
			\coordinate (A) at (0,0);
			\coordinate (B) at (0,1);
			\Rmorphism[A][$a$][$0$][][][][$g$];
			\ladder[B][][][$-lg$][$a+lkg$][$0$]
			\end{tikzpicture}
		}
	\end{array}
	:
	\begin{array}{c}
		\includeTikz{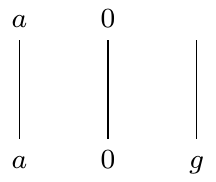}
		{
			\begin{tikzpicture}[xscale=.75]
			\coordinate (A) at (0,0);
			\Rladder[A][$a$][$0$][][$a$][$0$][$g$];
			\end{tikzpicture}
		}
	\end{array}
	\tilde{\to}
	\begin{array}{c}
		\includeTikz{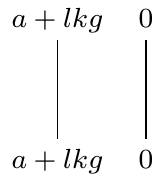}
		{
			\begin{tikzpicture}[xscale=.75]
			\coordinate (A) at (0,0);
			\ladder[A][$a+lkg$][$0$][][$a+lkg$][$0$];
			\end{tikzpicture}
		}
	\end{array}.\label{eqn:RF0rightiso2}
\end{align}
The right action is therefore $(a,0)\otimes g=(a+klg,0)$. Since the associators on $X_i$ are trivial, the morphisms commute and the resulting associator is trivial. These $\vvec{\ZZ{p}}$ actions and associator show that
\begin{align}
	X_k\otimes_{\vvec{\ZZ{p}}} X_l=X_{kl},
\end{align}
where $kl$ is taken modulo $p$. This is the first instance of an invertible bimodule we have seen. The bimodules $X_k$ therefore form a subgroup of the \emph{Brauer-Picard group} of $\vvec{\ZZ{p}}$.

\end{exmp}
\begin{center}\rule{.8\textwidth}{1pt}\end{center}

\begin{exmp}[$F_q\otimes_{\vvec{\ZZ{p}}} F_r$]
There is a single simple object in $\lad{F_q,F_r}{}$, with a $p$ dimensional algebra of morphisms
\begin{align}
\begin{array}{c}
\includeTikz{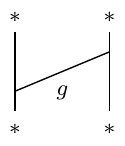}
{
	\begin{tikzpicture}[scale=0.8, every node/.style={scale=0.8}]
	\coordinate (A) at (0,0);
	\ladder[A][$*$][$*$][$g$][$*$][$*$];
	\end{tikzpicture}
}
\end{array}&,\\\text{with }
\begin{array}{c}
\includeTikz{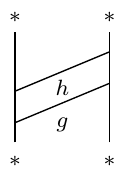}
{
	\begin{tikzpicture}[scale=0.8, every node/.style={scale=0.8}]
	\coordinate (A) at (0,0);\coordinate (B) at (0,.5);
	\ladder[A][$*$][$*$][$g$][][];
	\ladder[B][][][$h$][$*$][$*$];
	\end{tikzpicture}
}
\end{array}
&=\begin{array}{c}
\includeTikz{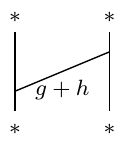}
{
	\begin{tikzpicture}[scale=0.8, every node/.style={scale=0.8}]
	\coordinate (A) at (0,0);
	\ladder[A][$*$][$*$][$g+h$][$*$][$*$];
	\end{tikzpicture}
}
\end{array}
\end{align}
which form the group algebra $\mathbb{C}\left[\ZZ{p}\right]$. The $p$ indecomposable idempotents $**\to **$ are given by
\begin{align}
I_{qk}&=\frac{1}{p}\sum_{g\in\ZZ{p}}
e^{\frac{2\pi i}{p}qkg}
\begin{array}{c}
\includeTikz{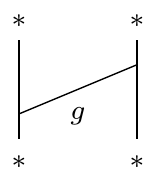}
{
	\begin{tikzpicture}
	\coordinate (A) at (0,0);
	\ladder[A][$*$][$*$][$g$][$*$][$*$];
	\end{tikzpicture}
}
\end{array}.
\end{align}
The $p$ simple objects of $\kar{\lad{F_q,F_r}{}}$ can be labelled by $k$. The morphism associated to the left $g$ action is
\begin{align}
\frac{1}{p}\sum_{x\in\ZZ{p}}e^{\frac{2\pi i}{p}qkx}
\begin{array}{c}
\includeTikz{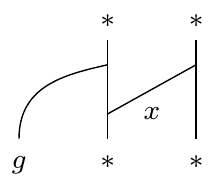}
{
	\begin{tikzpicture}[xscale=.75]
	\coordinate (A) at (0,0);
	\Lmorphism[A][$*$][$*$][$x$][$*$][$*$][$g$];
	\end{tikzpicture}
}
\end{array}
:
\frac{1}{p}\sum_{y\in\ZZ{p}}e^{\frac{2\pi i}{p}qky}
\begin{array}{c}
\includeTikz{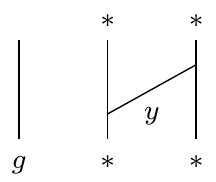}
{
	\begin{tikzpicture}[xscale=.75]
	\coordinate (A) at (0,0);
	\Lladder[A][$*$][$*$][$y$][$*$][$*$][$g$];
	\end{tikzpicture}
}
\end{array}
\tilde{\to}
\frac{1}{p}\sum_{y\in\ZZ{p}}e^{\frac{2\pi i}{p}q(k+g)y}
\begin{array}{c}
\includeTikz{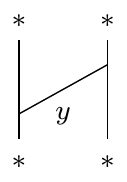}
{
	\begin{tikzpicture}[xscale=.75]
	\coordinate (A) at (0,0);
	\ladder[A][$*$][$*$][$y$][$*$][$*$];
	\end{tikzpicture}
}
\end{array}.\label{eqn:FFleftiso}
\end{align}
One can check that $I_{k+g}$ is absorbed from above
\begin{align}
\frac{1}{p^2}\sum_{x,y\in \ZZ{p}}
e^{\frac{2\pi i}{p}q(kx+(k+g)y)}
\begin{array}{c}
\includeTikz{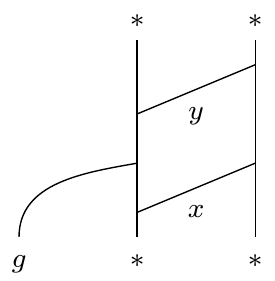}
{
	\begin{tikzpicture}
	\coordinate (A) at (0,0);
	\coordinate (B) at (0,1);
	\ladder[B][][][$y$][$*$][$*$];
	\Lmorphism[A][$*$][$*$][$x$][][][$g$];
	\end{tikzpicture}
}
\end{array}
&=
\frac{1}{p^2}\sum_{x,y\in \ZZ{p}}
e^{\frac{2\pi i}{p}q(kx+(k+g)y)}e^{-\frac{2\pi i}{p}qgy}
\begin{array}{c}
\includeTikz{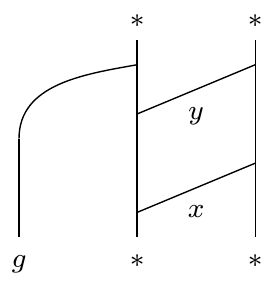}
{
	\begin{tikzpicture}
	\coordinate (A) at (0,0);
	\coordinate (B) at (0,-1);
	\Lladder[B][$*$][$*$][$x$][][][$g$];
	\Lmorphism[A][][][$y$][$*$][$*$][];
	\end{tikzpicture}
}
\end{array}\\
&=\frac{1}{p}\sum_{x\in \ZZ{p}}
e^{\frac{2\pi i}{p}qkx}
\begin{array}{c}
\includeTikz{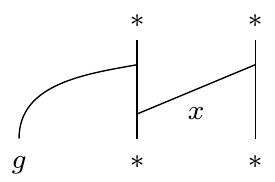}
{
	\begin{tikzpicture}
	\coordinate (A) at (0,0);
	\Lmorphism[A][$*$][$*$][$x$][$*$][$*$][$g$];
	\end{tikzpicture}
}
\end{array},
\end{align}
where the first equality uses the associator in Table~\ref{tab:zpdata}, where $y$ outgoing contributes $-y$. A similar computation on the right gives us 
\begin{align}
\frac{1}{p}\sum_{x\in\ZZ{p}}e^{\frac{2\pi i}{p}qkx}
\begin{array}{c}
\includeTikz{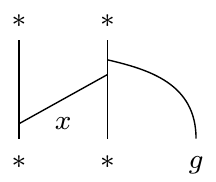}
{
	\begin{tikzpicture}[xscale=.75]
	\coordinate (A) at (0,0);
	\Rmorphism[A][$*$][$*$][$x$][$*$][$*$][$g$];
	\end{tikzpicture}
}
\end{array}
:
\frac{1}{p}\sum_{y\in\ZZ{p}}e^{\frac{2\pi i}{p}qky}
\begin{array}{c}
\includeTikz{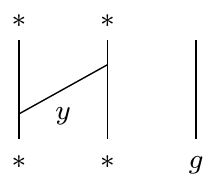}
{
	\begin{tikzpicture}[xscale=.75]
	\coordinate (A) at (0,0);
	\Rladder[A][$*$][$*$][$y$][$*$][$*$][$g$];
	\end{tikzpicture}
}
\end{array}
\tilde{\to}
\frac{1}{p}\sum_{y\in\ZZ{p}}e^{\frac{2\pi i}{p}q(k+q^{-1}rg)y}
\begin{array}{c}
\includeTikz{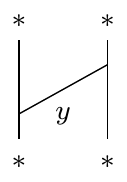}
{
	\begin{tikzpicture}[xscale=.75]
	\coordinate (A) at (0,0);
	\ladder[A][$*$][$*$][$y$][$*$][$*$];
	\end{tikzpicture}
}
\end{array}.\label{eqn:FFrightiso}
\end{align}
The left and right actions are therefore $g\otimes k=k+g$ and $k\otimes g=k+q^{-1}rg$. One can again check that the associator is trivial, giving us 
\begin{align}
	F_q\otimes_{\vvec{\ZZ{p}}} F_r=X_{q^{-1}r},
\end{align}
where $k^{-1}$ is the multiplicative inverse of $k$ taken modulo $p$. The bimodules $F_q$ are therefore invertible.

\end{exmp}

\begin{center}\rule{.8\textwidth}{1pt}\end{center}

\begin{exmp}[$F_q\otimes_{\vvec{\ZZ{p}}} X_l$]

The $p$ simple objects of $\lad{F_q,X_l}{}$ are labeled by $x\in\ZZ{p}$. Given any two objects $x$ and $y$, there is an isomorphism in $\lad{F_q,X_l}{}$, given by 
\begin{align}
\begin{array}{c}
\includeTikz{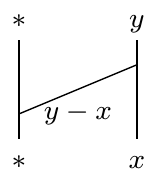}
{
	\begin{tikzpicture}
	\coordinate (A) at(0,0);
	\ladder[A][$*$][$x$][$y-x$][$*$][$y$];
	\end{tikzpicture}
}
\end{array}.
\end{align}
The single isomorphism class in $\kar{\lad{F_q,X_l}{}}$ can be represented using
\begin{align}
\begin{array}{c}
\includeTikz{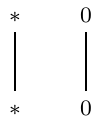}
{
	\begin{tikzpicture}[scale=0.6, every node/.style={scale=0.8}]
	\coordinate (A) at(0,0);
	\ladder[A][$*$][$0$][][$*$][$0$];
	\end{tikzpicture}
}
\end{array}.
\end{align}
Since there is a single object, we already know that $F_q\otimes_{\vvec{\ZZ{p}}}X_l=F_r$ for some $r$. To identify the value of $r$, we need to compute the associator. The left action is given by
\begin{align}
\begin{array}{c}
\includeTikz{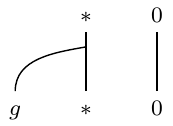}
{
	\begin{tikzpicture}[scale=0.6, every node/.style={scale=0.8}]
	\coordinate (A) at (0,0);
	\Lmorphism[A][$*$][$0$][][$*$][$0$][$g$];
	\end{tikzpicture}
}
\end{array}
:
\begin{array}{c}
\includeTikz{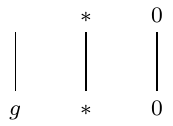}
{
	\begin{tikzpicture}[scale=0.6, every node/.style={scale=0.8}]
	\coordinate (A) at (0,0);
	\Lladder[A][$*$][$0$][][$*$][$0$][$g$];
	\end{tikzpicture}
}
\end{array}
\tilde{\to}
\begin{array}{c}
\includeTikz{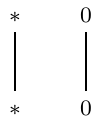}
{
	\begin{tikzpicture}[scale=0.6, every node/.style={scale=0.8}]
	\coordinate (A) at (0,0);
	\ladder[A][$*$][$0$][][$*$][$0$];
	\end{tikzpicture}
}
\end{array}.\label{eqn:FXleftiso}
\end{align}
On the right, we identify the isomorphism
\begin{align}
\begin{array}{c}
\includeTikz{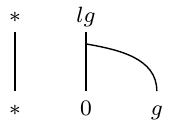}
{
	\begin{tikzpicture}[scale=0.6, every node/.style={scale=0.8}]
	\coordinate (A) at (0,0);
	\Rmorphism[A][$*$][$0$][][$*$][$lg$][$g$];
	\end{tikzpicture}
}
\end{array}
:
\begin{array}{c}
\includeTikz{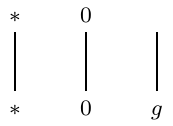}
{
	\begin{tikzpicture}[scale=0.6, every node/.style={scale=0.8}]
	\coordinate (A) at (0,0);
	\Rladder[A][$*$][$0$][][$*$][$0$][$g$];
	\end{tikzpicture}
}
\end{array}
\tilde{\to}
\begin{array}{c}
\includeTikz{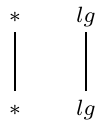}
{
	\begin{tikzpicture}[scale=0.6, every node/.style={scale=0.8}]
	\coordinate (A) at (0,0);
	\ladder[A][$*$][$lg$][][$*$][$lg$];
	\end{tikzpicture}
}
\end{array}.\label{eqn:FXrightiso}
\end{align}
The codomain is not the chosen representative of the isomorphism class. Instead, we need to append the morphism
\begin{align}
\begin{array}{c}
\includeTikz{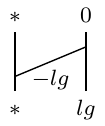}
{
	\begin{tikzpicture}[scale=0.6, every node/.style={scale=0.8}]
	\coordinate (A) at(0,0);
	\ladder[A][$*$][$lg$][$-lg$][$*$][$0$];
	\end{tikzpicture}
}
\end{array}
\end{align}
to correctly absorb the $g$. The complete isomorphism is given by
\begin{align}
\begin{array}{c}
\includeTikz{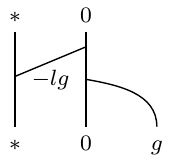}
{
	\begin{tikzpicture}[scale=0.6, every node/.style={scale=0.8}]
	\coordinate (A) at (0,0);
	\coordinate (B) at (0,1);
	\Rmorphism[A][$*$][$0$][][][][$g$];
	\ladder[B][][][$-lg$][$*$][$0$]
	\end{tikzpicture}
}
\end{array}
:
\begin{array}{c}
\includeTikz{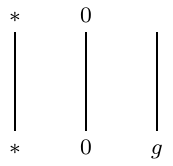}
{
	\begin{tikzpicture}[xscale=0.6, every node/.style={scale=0.8}]
	\coordinate (A) at (0,0);
	\Rladder[A][$*$][$0$][][$*$][$0$][$g$];
	\end{tikzpicture}
}
\end{array}
\tilde{\to}
\begin{array}{c}
\includeTikz{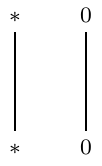}
{
	\begin{tikzpicture}[xscale=0.6, every node/.style={scale=0.8}]
	\coordinate (A) at (0,0);
	\ladder[A][$*$][$0$][][$*$][$0$];
	\end{tikzpicture}
}
\end{array}.\label{eqn:FXrightiso_full}
\end{align}
The associator $\omega$ can be identified by composing Equation~\ref{eqn:FXleftiso} and Equation~\ref{eqn:FXrightiso_full} in both orders
\begin{align}
\omega^{-1}(g,h)
\begin{array}{c}
\includeTikz{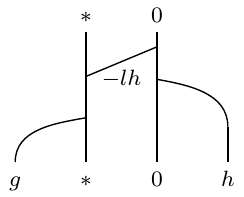}
{
	\begin{tikzpicture}[scale=0.6, every node/.style={scale=0.8}]
		\coordinate (A) at (0,-1);
		\coordinate (B) at (0,0);
		\coordinate (C) at (0,1);
		\Lmorphism[A][$*$][$0$][][][][$g$];
		\Rmorphism[B][][][][][][];
		\ladder[C][][][$-lh$][$*$][$0$];
		\draw(1.8,-.5)--(1.8,-1.1);
		\node[below,inline text] at (1.8,-1.1) {\footnotesize$h$};
	\end{tikzpicture}
}
\end{array}
&=
\begin{array}{c}
\includeTikz{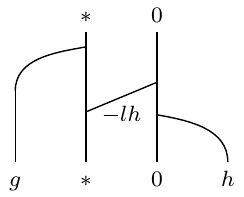}
{
	\begin{tikzpicture}[scale=0.6, every node/.style={scale=0.8}]
	\coordinate (A) at (0,-1);
	\coordinate (B) at (0,0);
	\coordinate (C) at (0,1);
	\Lmorphism[C][][][][$*$][$0$][];
	\Rmorphism[A][$*$][$0$][][][][$h$];
	\ladder[B][][][$-lh$][][];
	\draw(-1.8,.1)--(-1.8,-1.1);
	\node[below,inline text] at (-1.8,-1.1) {\footnotesize$g$};
	\end{tikzpicture}
}
\end{array}\\
&=e^{-\frac{2\pi i}{p}ql gh}
\begin{array}{c}
\includeTikz{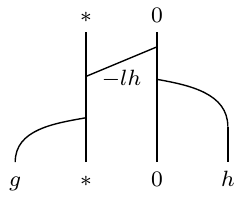}
{
	\begin{tikzpicture}[scale=0.6, every node/.style={scale=0.8}]
	\coordinate (A) at (0,-1);
	\coordinate (B) at (0,0);
	\coordinate (C) at (0,1);
	\Lmorphism[A][$*$][$0$][][][][$g$];
	\Rmorphism[B][][][][][][];
	\ladder[C][][][$-lh$][$*$][$0$];
	\draw(1.8,-.5)--(1.8,-1.1);
	\node[below,inline text] at (1.8,-1.1) {\footnotesize$h$};
	\end{tikzpicture}
}
\end{array},
\end{align}
where the final equality is obtained using the associator on the left strand. This gives us
\begin{align}
F_q\otimes_{\vvec{\ZZ{p}}} X_l=F_{ql}.
\end{align}

\end{exmp}

\begin{center}\rule{.8\textwidth}{1pt}\end{center}

\section{The topological phase $\vvec{\ZZ{p}}$}
Given a fusion category $\mathcal{C}$, the associated (nonchiral) topological phase is described by $Z(\mathcal{C})$, the \emph{Drinfeld center} of $\mathcal{C}$. The particles of the theory correspond to the simple objects of $Z(\mathcal{C})$, whilst the morphisms describe braiding processes. Given a modular category $Z(\mathcal{C})$, one can use the Levin-Wen procedure\cite{Levin2005} to construct a lattice model that realizes the topological phase in its low energy space.

Bimodule categories $\mathcal{C}\curvearrowright\mathcal{M}\curvearrowleft\mathcal{D}$ are associated to domain walls between the topological phases $Z(\mathcal{C})$ and $Z(\mathcal{D})$. In this section, we discuss a physical interpretation of the bimodule categories in Table~\ref{tab:zpdata}. We also show how the wall fusion rules, including the multiplicities, in Table~\ref{tab:zptable} can be obtained from the physical theory. 


The lattice model associated to $Z(\vvec{\ZZ{p}})$ has $p^2$ particle-like excitations, usually labelled $e^am^b$, with $a,b\in\ZZ{p}$. The fusion rules are $e^a\times e^b=e^{a+b}$ and $m^a\times m^b=m^{a+b}$ and the braiding is defined by $e^a m^b=e^{\frac{2\pi i}{p}ab}m^be^a$. On a torus, this model has $p^2$ degenerate ground states which can be used to encode a pair of $p$-dimensional logical spaces (referred to a qu$p$its). The logical Pauli gates are implemented by moving a particle around one of the cycles of the torus, so the encoded information is protected against local noise.

There are two possible boundary conditions that can be imposed on the model, commonly referred to as rough and smooth due to their lattice implementation\cite{MR2942952}. At a rough boundary, a particle $e^a$ can be `condensed' into the vacuum. In other words, fusing an $e$ particle into the boundary leaves the boundary unaffected. The smooth boundaries behave similarly for the $m$ sector. A single qu$p$it can be encoded into a strip of the phase with smooth (or rough) boundaries on both sides. By creating a hole or puncture to vacuum in the middle of the code, additional qu$p$its can be encoded\cite{0610082}. By moving these punctures around each other, logical gates can be applied to the encoded information\cite{0610082,Bombin2007a,Brown2016}.  

The possible invertible domain walls for this theory are well understood\cite{Bridgeman2017}. These walls form the group $\dih{p-1}$, the dihedral group of order $2p-2$. By sweeping these walls across the lattice, additional logical gates can be applied to the encoded information. 

By allowing interfaces between distinct domain walls, such as corners between rough and smooth boundaries, or terminations of invertible domain walls, additional qu$p$its can be encoded and manipulated\cite{Bombin2010,1606.07116,Brown2016}. A comprehensive study of the possible \emph{binary interface defects} is the subject of the companion paper \onlinecite{DefectsPaper}.

The group of invertible domain walls, along with the boundaries can be used to realize all bimodules. A domain wall can be built from a pair of boundaries by simply stacking them, for example a smooth and rough boundary
\begin{align}
	\begin{array}{c}
		\includeTikz{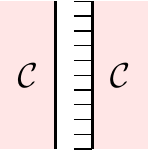}
		{
			\begin{tikzpicture}[scale=.75]
				\fill[red!10](-1,-1) rectangle (-.25,1);
				\draw[thick](-.25,-1)--(-.25,1);
				\fill[red!10](1,-1) rectangle (.25,1);
				\draw[thick](.25,-1)--(.25,1);
				\foreach \x in {0,1,...,10}
				{
					\draw[shift={(0,1-.2*\x)}](0,0)--(.25,0);
				}
				\node at (-.625,0) {$\mathcal{C}$};\node at (.625,0) {$\mathcal{C}$};
			\end{tikzpicture}
		}
	\end{array}.
\end{align}

\begin{table}
\begin{tabular}{!{\vrule width 1pt}c!{\vrule width 1pt}c|c!{\vrule width 1pt}}
	\toprule[1pt]
	\rowcolor[gray]{.9}[\tabcolsep]Bimodule label & Domain wall & Action on particles\\
	\toprule[1pt]
	$T$&$\begin{array}{c}\includeTikz{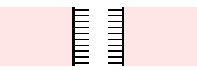}{
		\begin{tikzpicture}[yscale=.3]
		\draw[white](0,-1.1)--(0,1.21);
		\fill[red!10](-1,-1) rectangle (-.25,1);
		\fill[red!10](1,-1) rectangle (.25,1);
		\draw[thick](-.25,-1)--(-.25,1);
		\draw[thick](.25,-1)--(.25,1);
		\foreach \x in {0,...,9}{\draw (.1,-.9+.2*\x)--(.25,-.9+.2*\x);};
		\foreach \x in {0,...,9}{\draw (-.1,-.9+.2*\x)--(-.25,-.9+.2*\x);};
		\end{tikzpicture}}
	\end{array}$&Condenses $e$ on both sides\\
	$L$&$\begin{array}{c}\includeTikz{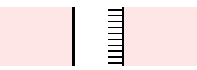}{\begin{tikzpicture}[yscale=.3]
		\draw[white](0,-1.1)--(0,1.21);
		\fill[red!10](-1,-1) rectangle (-.25,1);
		\fill[red!10](1,-1) rectangle (.25,1);
		\draw[thick](-.25,-1)--(-.25,1);
		\draw[thick](.25,-1)--(.25,1);
		\foreach \x in {0,...,9}{\draw (.1,-.9+.2*\x)--(.25,-.9+.2*\x);};
		\end{tikzpicture}}\end{array}$&Condenses $m$ on left and $e$ on right\\
	$R$&$\begin{array}{c}\includeTikz{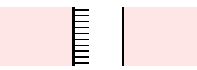}{\begin{tikzpicture}[yscale=.3]
		\draw[white](0,-1.1)--(0,1.21);
		\fill[red!10](-1,-1) rectangle (-.25,1);
		\fill[red!10](1,-1) rectangle (.25,1);
		\draw[thick](-.25,-1)--(-.25,1);
		\draw[thick](.25,-1)--(.25,1);
		\foreach \x in {0,...,9}{\draw (-.1,-.9+.2*\x)--(-.25,-.9+.2*\x);};
		\end{tikzpicture}}\end{array}$&Condenses $e$ on left and $m$ on right\\
	$F_0$&$\begin{array}{c}\includeTikz{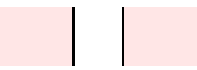}{\begin{tikzpicture}[yscale=.3]
		\draw[white](0,-1.1)--(0,1.21);
		\fill[red!10](-1,-1) rectangle (-.25,1);
		\fill[red!10](1,-1) rectangle (.25,1);
		\draw[thick](-.25,-1)--(-.25,1);
		\draw[thick](.25,-1)--(.25,1);
		\end{tikzpicture}}\end{array}$&Condenses $m$ on both sides\\
	\toprule[1pt]
	$X_k$&$\begin{array}{c}\includeTikz{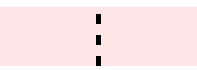}{\begin{tikzpicture}[yscale=.3]
		\draw[white](0,-1.1)--(0,1.21);
		\fill[red!10](-1,-1) rectangle (1,1);
		\draw[ultra thick,dashed] (0,-1)--(0,1);
		\end{tikzpicture}}\end{array}$&$X_k:e^am^b\mapsto e^{ka}m^{k^{-1} b}$ (moving left to right), where $k^{-1}$ is taken multiplicatively modulo $p$\\
	$F_{q}=F_1 X_q$&$\begin{array}{c}\includeTikz{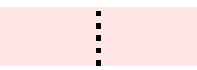}{\begin{tikzpicture}[yscale=.3]
		\draw[white](0,-1.1)--(0,1.21);
		\fill[red!10](-1,-1) rectangle (1,1);
		\draw[ultra thick,dotted] (0,-1)--(0,1);\end{tikzpicture}}\end{array}$&
	$F_1:e^am^b\mapsto e^{b}m^{a}$\\
	\toprule[1pt]
\end{tabular}
	\caption{Domain walls on the lattice corresponding to bimodules}\label{tab:zpinterp}
\end{table}

Since there are 2 types of boundary, there are 4 non-invertible domain walls in $\vvec{\ZZ{p}}$. In Table~\ref{tab:zpinterp}, we associate a domain wall to each bimodule. The fusion of bimodules can be understood via this association. For example, 
\begin{align}
	T\otimes_{\vvec{\ZZ{p}}} T\to
	\begin{array}{c}
		\includeTikz{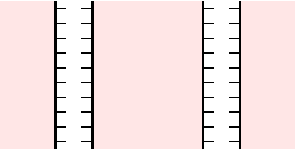}
		{
			\begin{tikzpicture}[scale=.75]
				\fill[red!10] (-2,-1)rectangle(-1.25,1);
				\draw[thick] (-1.25,-1)--(-1.25,1);
				\fill[red!10] (-.75,-1)rectangle(.75,1);
				\draw[thick,fill=none] (.75,-1)--(.75,1) (-.75,1)--(-.75,-1);
				\fill[red!10] (2,-1)rectangle(1.25,1);
				\draw[thick] (1.25,-1)--(1.25,1);
				\foreach \x in {0,...,9}
				{
					\draw[shift={(0,.1+.2*\x)}] (-1.25,-1)--(-1.1,-1);
					\draw[shift={(0,.1+.2*\x)}] (-.75,-1)--(-.9,-1);
					\draw[shift={(0,.1+.2*\x)}] (1.25,-1)--(1.1,-1);
					\draw[shift={(0,.1+.2*\x)}] (.75,-1)--(.9,-1);
				}
			\end{tikzpicture}
		}
	\end{array}
	&\propto
	\begin{array}{c}
		\includeTikz{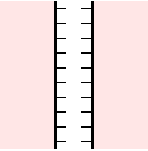}
		{
			\begin{tikzpicture}[scale=.75]
			\fill[red!10] (-1,-1)rectangle(-.25,1);
			\draw[thick] (-.25,-1)--(-.25,1);
			\fill[red!10] (1,-1)rectangle(.25,1);
			\draw[thick] (.25,-1)--(.25,1);
			\foreach \x in {0,...,9}
				{
					\draw[shift={(0,.1+.2*\x)}] (-.25,-1)--(-.1,-1);
					\draw[shift={(0,.1+.2*\x)}] (.25,-1)--(.1,-1);
				}
			\end{tikzpicture}
		}
	\end{array},
\end{align}
by discarding the isolated strip. The theory $Z(\vvec{\ZZ{p}})$ on a strip or tube with both rough (or both smooth) boundaries supports a $p$ dimensional vector space. These $p$ distinct states for the central strip give the coefficient $p$. On the other hand, if the opposite boundaries are not the same, the state is unique and the coefficient is 1 
\begin{align}
	T\otimes_{\vvec{\ZZ{p}}} L\to
	\begin{array}{c}
		\includeTikz{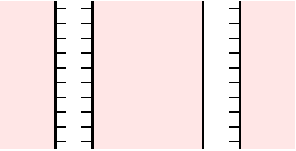}
		{
			\begin{tikzpicture}[scale=.75]
			\fill[red!10] (-2,-1)rectangle(-1.25,1);
			\draw[thick] (-1.25,-1)--(-1.25,1);
			\fill[red!10] (-.75,-1)rectangle(.75,1);
			\draw[thick,fill=none] (.75,-1)--(.75,1) (-.75,1)--(-.75,-1);
			\fill[red!10] (2,-1)rectangle(1.25,1);
			\draw[thick] (1.25,-1)--(1.25,1);
			\foreach \x in {0,...,9}
			{
				\draw[shift={(0,.1+.2*\x)}] (-1.25,-1)--(-1.1,-1);
				\draw[shift={(0,.1+.2*\x)}] (-.75,-1)--(-.9,-1);
				\draw[shift={(0,.1+.2*\x)}] (1.25,-1)--(1.1,-1);
			}
			\end{tikzpicture}
		}
	\end{array}
	&=
	\begin{array}{c}
		\includeTikz{TxTpic2}
		{
			\begin{tikzpicture}[scale=.75]
			\fill[red!10] (-1,-1)rectangle(-.25,1);
			\draw[thick] (-.25,-1)--(-.25,1);
			\fill[red!10] (1,-1)rectangle(.25,1);
			\draw[thick] (.25,-1)--(.25,1);
			\foreach \x in {0,...,9}
			{
				\draw[shift={(0,.1+.2*\x)}] (-.25,-1)--(-.1,-1);
				\draw[shift={(0,.1+.2*\x)}] (.25,-1)--(.1,-1);
			}
			\end{tikzpicture}
		}
	\end{array}.
\end{align}

The action of the invertible domain walls corresponding to $F_q$ swap the $e$ and $m$ sectors of the theory. The effect on the boundaries is to switch rough and smooth. From this, we obtain the same multiplication Table~\ref{tab:zptable}.

It follows from our computations that every domain wall between the $\protect\vvec{\ZZ{p}}$ theory and itself are either invertible, or have a strip of vacuum separating the two bulks. If $p \neq q$, then all the domain walls between the $\protect\vvec{\ZZ{p}}$ and $\protect\vvec{\ZZ{q}}$ theories have a strip of vacuum separating the two bulks. This is consistent with the computation using ladder string diagrams.

%

\section{Conclusions} \label{S:conc}

In this work, we have studied domain walls between topological phases. Using ladder string diagrams, we have shown how the relative tensor product of bimodule categories can be computed. We have verified that the relative tensor product corresponds to the stacking of domain walls in the $\vvec{\ZZ{p}}$ case. Physically, this fusion category corresponds to the quantum double model of $\ZZ{p}$ and is important for both quantum memories and quantum information processing tasks.

Although we have only considered $\vvec{\ZZ{p}}$, the framework outlined here is not restricted to this class of fusion categories. The color code\cite{ColorCode,Yoshida2015a} ($\vvec{\ZZ{2}\times\ZZ{2}}$) is a important quantum code due to the number of fault-tolerant gates that can be implemented. This theory has nontrivial domain walls to both the $\vvec{\ZZ{2}}$ model (toric code) and $\vvec{\ZZ{4}}$ model. Although not fully invertible, we expect domain walls which do not have a `vacuum gap' (unlike $T$ for example) between these theories. Is would be interested to investigate these walls, and see if they extend the computational power of these models. Moreover, we expect that the developed tools will be useful for computing the Brauer-Picard ring for $\bpr{\vvec{G}}$ when $G$ is not abelian.

The domain walls considered in this work did not terminate, rather they went on to infinity or had periodic boundaries. The theory of defects of invertible domain walls\cite{Barkeshli2013,Barkeshli2014,SETPaper} has already found application in quantum information tasks\cite{Bombin2010,Brown2016}. It is therefore important to ask how more general defects can be utilized. The categorical framework used in this paper is useful for understanding these defects\cite{DefectsPaper}. Indeed, defects are modeled using functors between bimodules.

\acknowledgments
This work is supported by the Australian Research Council (ARC) via the Centre of Excellence in Engineered Quantum Systems (EQuS) project number CE170100009, Discovery Project ``Subfactors and symmetries`` DP140100732 and Discovery Project ``Low dimensional categories`` DP16010347. We thank Andrew Doherty, Steve Flammia, Dominic Williamson, Cain Edie-Michell, Scott Morrison and Christoph Schweigert. This work would not have been possible without their input and feedback. We thank Paul Wedrich and Christopher Chubb for feedback on the draft manuscript. We thank Noah Snyder for explaining how to prove proposition \ref{eno_remark_3.9}.

\vspace{5mm}

\bibliographystyle{apsrev_jacob}
\bibliography{refs}

\end{document}